\documentclass[review,3p]{elsarticle}
\usepackage{subfigure}
\usepackage{comment}
\usepackage{bm}
\usepackage{graphicx}
\usepackage{wrapfig}
\usepackage{booktabs}
\usepackage[colorlinks,linkcolor=blue]{hyperref}
\UseRawInputEncoding
\biboptions{numbers,sort&compress}
 \usepackage{mathrsfs}
\usepackage{makecell}
\usepackage{ulem}
\usepackage{mathrsfs}
\usepackage{appendix} 
\usepackage{setspace}
\usepackage{algorithmic} 
\usepackage[ruled]{algorithm2e}
\usepackage{multirow} 
\usepackage{xcolor}

\usepackage{amssymb}
\usepackage{amsthm}
\usepackage{amsmath}
\newcommand{\mt}{\mathbb{T}}

\journal{Computer methods in applied mechanics and engineering}

\begin{document}

\begin{frontmatter}

\title{Adaptive Isogeometric Topology Optimization of Shell Structures based on PHT-splines }

\author[vt]{Zepeng Wen}

\author[rvt]{Qiong Pan}

\author[rvt]{Xiaoya Zhai\corref{cor1}}
\ead{xiaoyazhai@ustc.edu.cn}

\author[vt]{Hongmei Kang\corref{cor1}}
\ead{khm@suda.edu.cn}
\cortext[cor1]{Corresponding author}

\author[rvt]{Falai Chen}%\corref{cor1}}
%\ead{chenfl@ustc.edu.cn}

\address[vt]{School of Mathematical Sciences, Soochow University, No.1 Road Shizi, Suzhou, 215006, P. R. China}
\address[rvt]{School of Mathematical Sciences, University of Science and Technology of China, No.96 Road jinzhai, Hefei, Anhui, 230026, P. R. China}

\begin{abstract}
This paper proposes an Adaptive Isogeometric Topology Optimization framework for shell structures based on PHT-splines (PHT-AITO). In this framework, the design domain, displacement, and density are represented by PHT-splines. 
Leveraging the local refinement capability of PHT-splines, mesh elements defining the density function are adaptively refined to achieve a suitable resolution at the interface between solid and void regions. 
This addresses the issue of excessive degrees of freedom resulting from global refinement. 
The refinement of the mesh elements is driven by their density. During the optimization of the density on a refined mesh, the initial value of the density is inherited from the optimization results on the previous mesh to accelerate the iteration process and maintain the stability of the optimized structure. Numerical experiments on various shell structures have verified the effectiveness of PHT-AITO. Compared with isogeometric topology optimization based on tensor-product splines, PHT-AITO can significantly reduce the degrees of freedom in the optimization problem, thereby improving computational efficiency.
\end{abstract}

\begin{keyword} 
Adaptive isogeometric analysis \sep PHT-splines \sep isogeometric topology optimization \sep shell structures
\end{keyword}
\end{frontmatter}

\section{Introduction}
\label{sec:Intro}

Shell structures are curved surfaces with thickness that have a wide range of applications in architecture and various industries. These structures are characterized by their lightweight design, stability, and aesthetic value~\cite{MICHAELROTTER19983}. In architecture, shell structures offer an optimal means of establishing open and flexible spaces. The incorporation of curved configurations and diverse material allocations frequently yields architecturally captivating and visually compelling concepts. In industrial applications, shell structures can attain exceptional strength-to-mass ratios while inherently manifesting stability rooted in their inherent material allocations. Consequently, optimizing the material distribution within the shell structures is a pivotal domain for scholarly exploration.

Topology optimization~\cite{bendsoe2003topology} is a computational design technique used to determine the optimal material distribution for achieving desired performance. It serves as a bridge between practical engineering scenarios and creative art design. The design and analysis of shell structures~\cite{farshad2013design} have been developed based on elastic shell theory~\cite{naghdi1962foundations} and nonlinear elastic shell theory~\cite{budiansky1968notes}. Ahmad et al.~\cite{ahmad1970analysis} first transferred the topology optimization of shell structures into the optimization of their mid-surfaces using curved elements. Subsequent important generalizations include optimization of the shell geometry and thickness by Hinton et al.~\cite{hinton1992finite} and Maute et al.~\cite{maute1997adaptive}.

Many well-developed topology optimization methods are applied to the topology optimization of shell structures. In \cite{ansola2002integrated,hassani2013simultaneous, Ye2019Topology}, the shape of the shell structure and material distribution are optimized simultaneously, particularly, the topology optimization stage of shell structures adopted the homogenization method~\cite{soto1993modelling}, the Solid Isotropic Material with Penalization Method (SIMP)~\cite{SIMP} and the Level Set Method (LSM)~\cite{levelset}, respectively. In~\cite{Townsend2019ALS}, the buckling behavior of the shell structure is taken as the optimization goal, and LSM is employed to optimize the shell structure to find the material layout with good buckling performance. The Moving Morphable Component (MMC) method~\cite{MMC} combined with conformal mapping is utilized for shell topology optimization in~\cite{Huo2022topology} to reduce the number of design variables. In~\cite{jiang2023explicit}, MMC is employed for the shape and topology optimization of shell structures, in which both the shell structure and a group of components are represented using Non-Uniform Rational B-spline (NURBS).

In general topology optimization, unknown structural responses are numerically solved by the Finite Element Method (FEM). However, the FEM approach has its limitations, as it necessitates the discretization of structural geometry into meshes and characterizes structural responses by low-order continuity ($C^0$). These limitations are known to be prone to causing numerical issues. Isogeometric Topology Optimization (ITO) combines the benefits of isogeometric analysis (IGA)~\cite{IGA,cottrell2009isogeometric} and topology optimization to offer significant advantages in the optimization of structural systems. By integrating CAD and analysis seamlessly, isogeometric topology optimization ensures accurate geometric representation without relying on discrete approximations. Moreover, it facilitates the construction of higher-order basis functions for analysis. 

The pioneering work of ITO might go back to \cite{seo2010shape,seo2010isogeometric}, where the topology optimization is conducted within the IGA framework using trimming techniques. The density-based variant of ITO was initially explored in \cite{hassani2012isogeometrical} and further improved in  \cite{gao2019isogeometric}, incorporating an enhanced density distribution function. Another approach combining LSM and IGA is proposed in \cite{shojaee2012composition, wang2016isogeometric} to conduct the structure optimization problem. In the MMC/MMV (Moving Morphable Void)-based ITO \cite{hou2017explicit, zhang2020explicit}, the parameters used for representing geometries of structural components (a subset of the design domain) are taken as design variables, consequently, this approach yields explicit structural representations. For shell structures, the ITO methods based on trimmed splines \cite{kang2016isogeometric} and MMV ~\cite{zhang2020explicit, zhang2020stress} are utilized for shell topology optimization. For more references to ITO, we recommend readers refer to \cite{gao2020comprehensive, itobook}.

The high-resolution design domain provides more extensive material layout options but requires time-consuming computations. To alleviate this, a strategy of employing finer elements along structure boundaries and coarser elements within void and solid regions effectively reduces computational complexity. Locally refinable splines \cite{lrsplines} have emerged as an amendment of NURBS. NURBS has limitations in representing complex geometric models due to its tensor product topological restrictions. The incorporation of locally refinable splines into isogeometric topology optimization, referred to as Adaptive Isogeometric Topology Optimization(AITO), enables the handling of design problems involving complicated structures while reducing the degrees of freedom required for high-resolution density field representation. More recently, various local refinable spline types, including THB-splines \cite{THBTO}, T-splines \cite{TTO}, and PHT-splines \cite{gupta2022adaptive,karuthedath2023continuous}, have been integrated into topology optimization, contributing to an adaptive ITO framework.

Among various locally refinable splines \cite{HBiga,Tsplines,AST,THB,LR,PHT}, PHT-splines (Polynomial splines on Hierarchical T-meshes) \cite{PHT} stand out for their concise and efficient local refinement algorithm and a set of analysis suitable basis functions. Additionally, PHT-splines provide a natural B\'ezier expression, eliminating the need for B\'ezier extraction in IGA. These advantages have led to successful application of PHT-splines in linear elastic problems \cite{PHTIGA,PHTthinshell,PHTelastic} and vibration problems \cite{PHTVibrations,PHTKL} within the IGA framework.

In this paper, we introduce an adaptive isogeometric topology optimization framework for shell structures based on PHT-splines. Within this framework, the shell structure and density function are represented by PHT-splines. The compliance of the shell structure is minimized under a volume constraint, where the control coefficients of the density function are used as design variables.
The optimization process begins with a coarse mesh and adaptively refines the mesh based on the density value at the element center. More precisely, elements undergo refinement when their density falls within the range of $[0.1,0.9]$. This refinement aims to achieve a high-resolution density representation, especially in the interface between solid and void regions. The shell optimization problem is resolved at each refined mesh level, with initial design variable values inherited from the optimized results on the previous mesh. The sensitivity filter, defined by the weighted sum of density sensitivity within the support of the relevant design variables, is used in optimization to update design variables. The requirement for a predefined filter radius in traditional methods has been eliminated.
Several representative shell structures are provided with different loads and boundary conditions to validate the effectiveness and robustness of the proposed method.

The contributions of this paper are concluded as follows:
\begin{itemize}
	\item We propose an Adaptive Isogeometric Topology Optimization framework for shell structures based on PHT-splines, called PHT-AITO. The density distribution is expressed as a continuous PHT-spline function with its control coefficients serving as design variables. PHT-AITO employs an adaptive mesh for the optimization of density distribution, featuring finer elements at the interface between solid and void material, and coarser elements within the solid and void regions. 
	\item For PHT-AITO, we propose the adaptive smoothing filter, in which the radius of the filter can be automatically and adaptively determined without specifying parameters in advance. Additionally, the density inheritance strategy is utilized in the optimization iterative process to ensure the reuse of the optimized results from different mesh size, aiming to improve the iterative optimization efficiency and the topological stability of the shell structure.
	\item Leveraging the local refinement capability of PHT-splines, PHT-AITO addresses the problem of excessive degrees of freedom caused by high-resolution tensor product meshes. In other words, this approach yields a high-quality density distribution with reduced degrees of freedom.
\end{itemize}

The paper is organized as follows. Section~\ref{sec:pre} provides the preliminary knowledge of PHT-splines and shell structures.
In Section~\ref{sec:method}, the adaptive isogeometric topology optimization based on PHT-splines is presented in detail.
Section~\ref{sec:res} demonstrates numerical experiments on various shell structures using the proposed PHT-AITO framework. Section~\ref{sec:con} concludes the paper and discusses the future potential work.

\section{Preliminaries}
\label{sec:pre}
In this section, we present preliminary knowledge related to PHT-splines and shell structures. These fundamentals are essential to build our proposed adaptive isogeometric analysis framework for shell structures.

\subsection{PHT-splines \label{sec.pht}}
PHT-splines (Polynomial splines over Hierarchical T-meshes), as a kind of locally refinable splines, provide an efficient local refinement algorithm and inherit numerous advantageous properties from T-splines. As a result, PHT-splines are widely applied in adaptive geometric modeling and isogeometric analysis. In this paper, PHT-splines are used for adaptive topology optimization.

A {\em T-mesh} is a rectangular grid that allows T-junctions. It is assumed that the endpoints of each grid line in the T-mesh must be on two other grid lines and each cell in the mesh must be a rectangle. The corners of all the rectangles constitute the \textit{vertices} of the T-mesh. The rectangle with no other edges and vertices inside is called an \textit{element} of the T-mesh. Vertices that are inside the mesh are called \textit{interior vertices}, otherwise, they are called \textit{boundary vertices}. Interior vertices can be \textit{crossing vertices} or \textit{T vertices}, depending on whether they cross or terminate at any grid lines. A \textit{hierarchical T-mesh} is a special type of T-mesh that has a natural level structure. It is defined in a recursive fashion. One generally starts from a tensor-product mesh (level $0$). From level $k$ to level $k + 1$, one subdivides a cell at level $k$ into four sub-elements which are elements at level $k +1$. For simplicity, we consider dyadic refinement.

Given a T-mesh $\mt$, $\mathcal{F}$ denotes all the elements in $\mt$ and $\Omega$ denotes the region occupied by $\mathcal{F}$. The polynomial spline space over $\mt$ is defined as
\begin{align}
	\mathcal{S}(m,n,\alpha, \beta,
	\mt) := & \Big\{ s(x,y) \in C^{\alpha,\beta}(\Omega)\Big|
	s(x,y)|_\phi \in \mathbb{P}_{mn} \mbox{ for any } \phi \in \mathcal{F}\Big\},
\end{align} 
where $\mathbb{P}_{mn}$ is the space of all the polynomials with bi-degree $(m,n)$, and $C^{\alpha, \beta}(\Omega)$ is the space consisting of all the bivariate functions which are continuous in $\Omega$ with order $\alpha$ in the $x$-direction and with order $\beta$ in the $y$-direction.

PHT-splines are exactly the $\mathcal{S}(3,3,1,1,\mt)$ with $\mt$ being a hierarchical T-mesh, which has a concise dimension formula
\begin{equation}
	\label{E.dimen}
	\mathcal{S}(3,3,1,1, \mt)=4(V^b+V^+),
\end{equation}
where $V^b$ and $V^+$ represent the number of boundary vertices and interior crossing vertices in $\mt$ respectively. A boundary vertex or an interior crossing vertex is called a \emph{basis vertex}. Each basis vertex is associated with four basis functions.

PHT-spline basis functions are constructed level-by-level \cite{bib:PHT}. For the initial level (denoted as $\mt_0$), the standard bicubic $C^1$ continuous tensor-product B-splines are used as basis functions. A hierarchical T-mesh at level $k$ is denoted as $\mt_k$. Suppose the basis functions $\{b^k_j\}$, $j=1,...,d_k$, on $\mt_k$ have been
constructed, then the basis functions on $\mt_{k+1}$ are constructed by two steps:
1) truncating the basis functions $\{b^k_j\}_{j=1}^{d_k}$ on $\mt_k$; 2) constructing bicubic $C^1$ continuous B-splines basis functions associated with the new basis vertices in $\mt_{k+1}$.

\begin{figure*}[htbp]
	\centering 
	\subfigure[The four basis functions associated with the basis vertex marked by a black dot on the initial mesh]{
		\includegraphics[width = .3\textwidth]{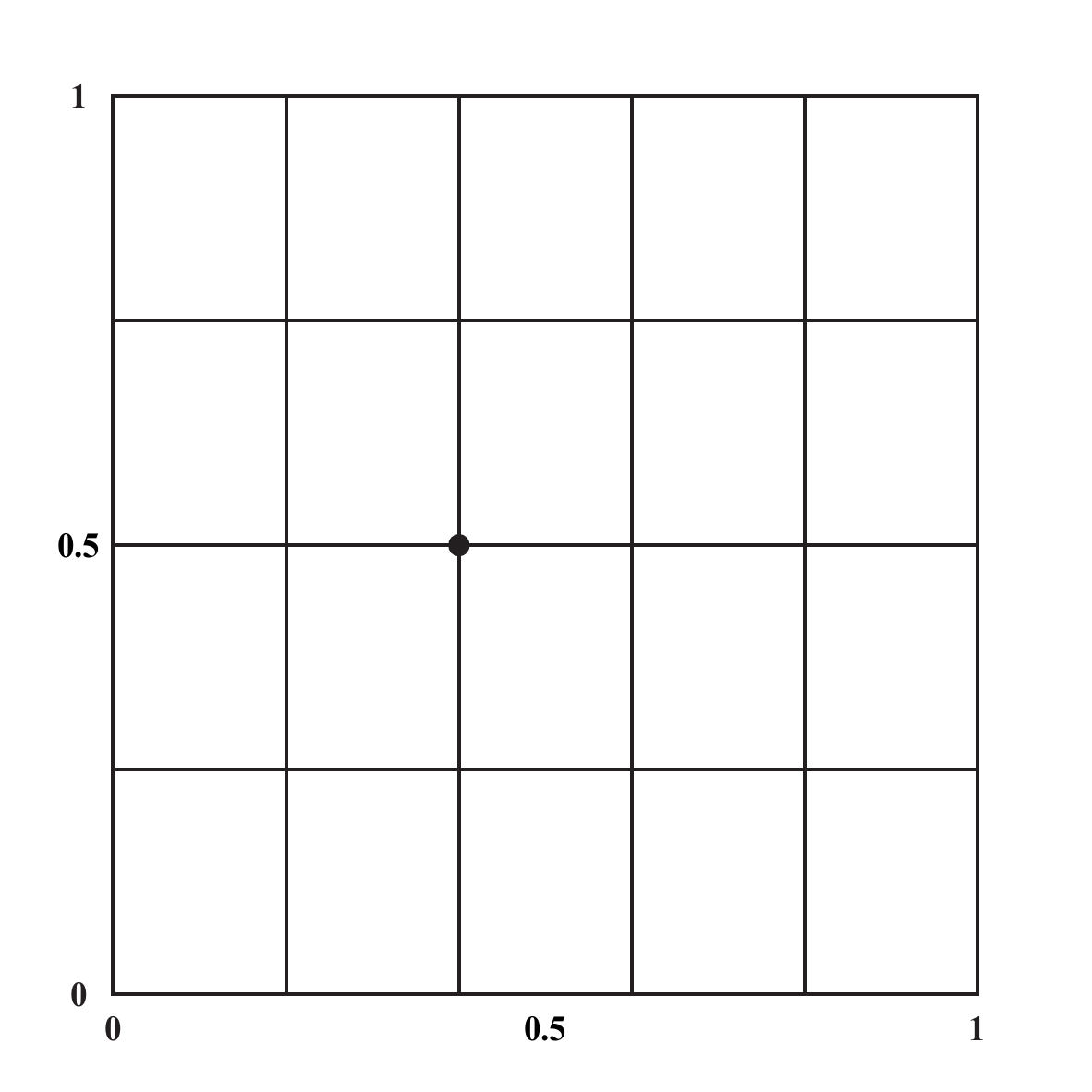}
		\includegraphics[width = .6\textwidth]{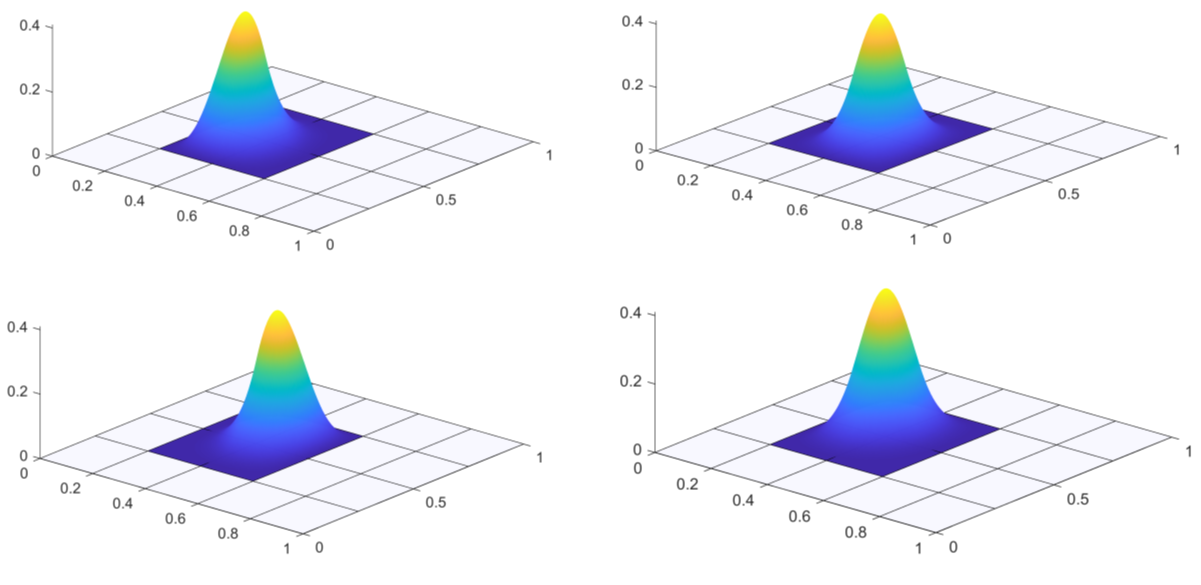}} 
	\subfigure[The four basis functions associated with the marked basis vertex on the refined mesh ]{
		\includegraphics[width = .3\textwidth]{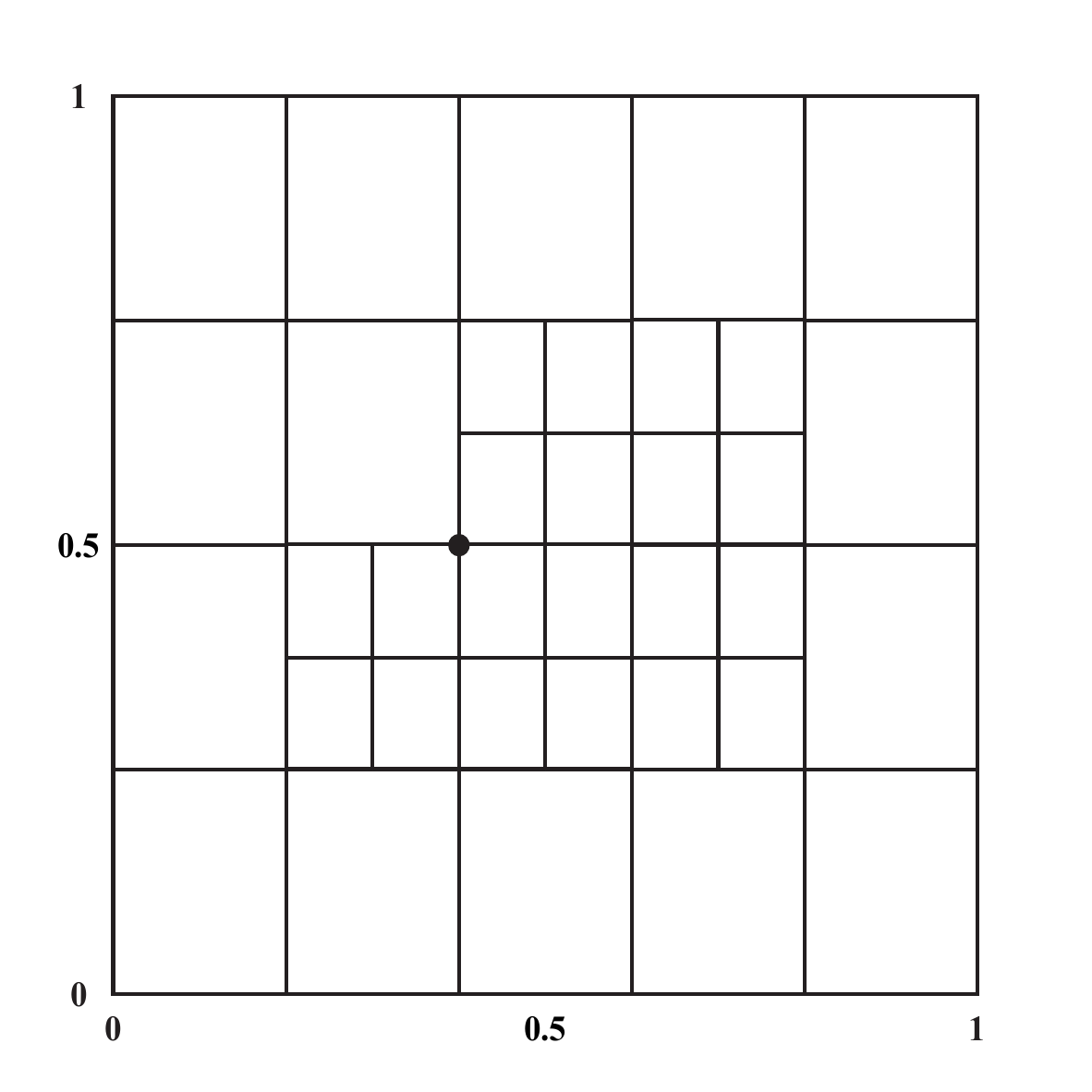}
		\includegraphics[width = .6\textwidth]{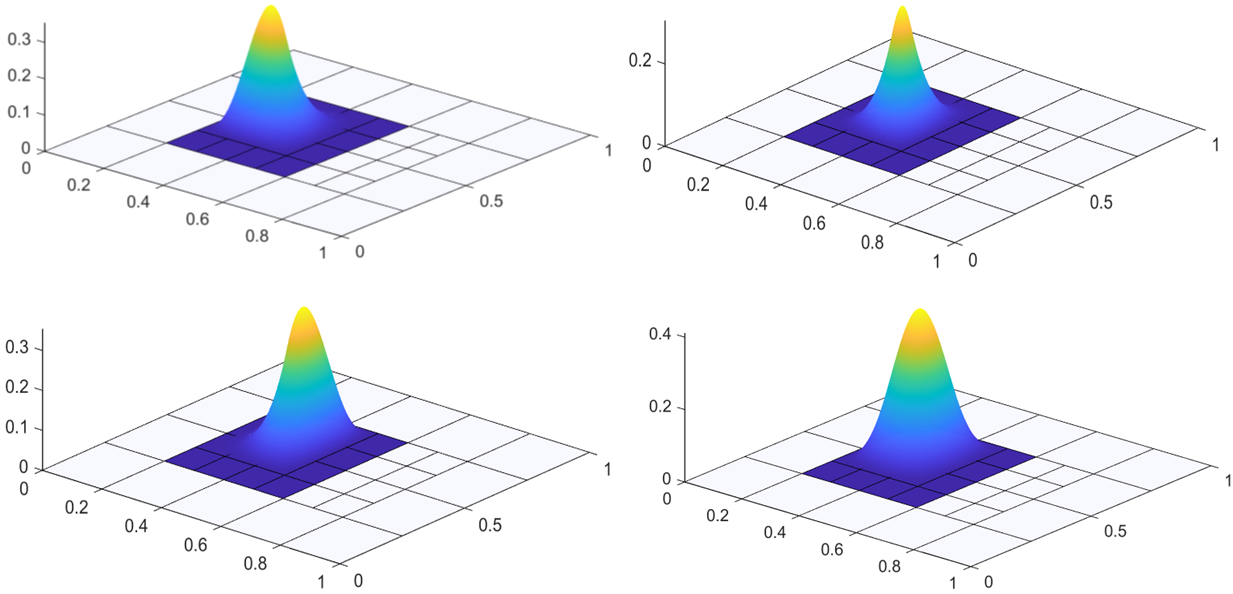}} 
	\subfigure[B\'{e}zier ordinates over the support]{
		\includegraphics[width = .3\textwidth]{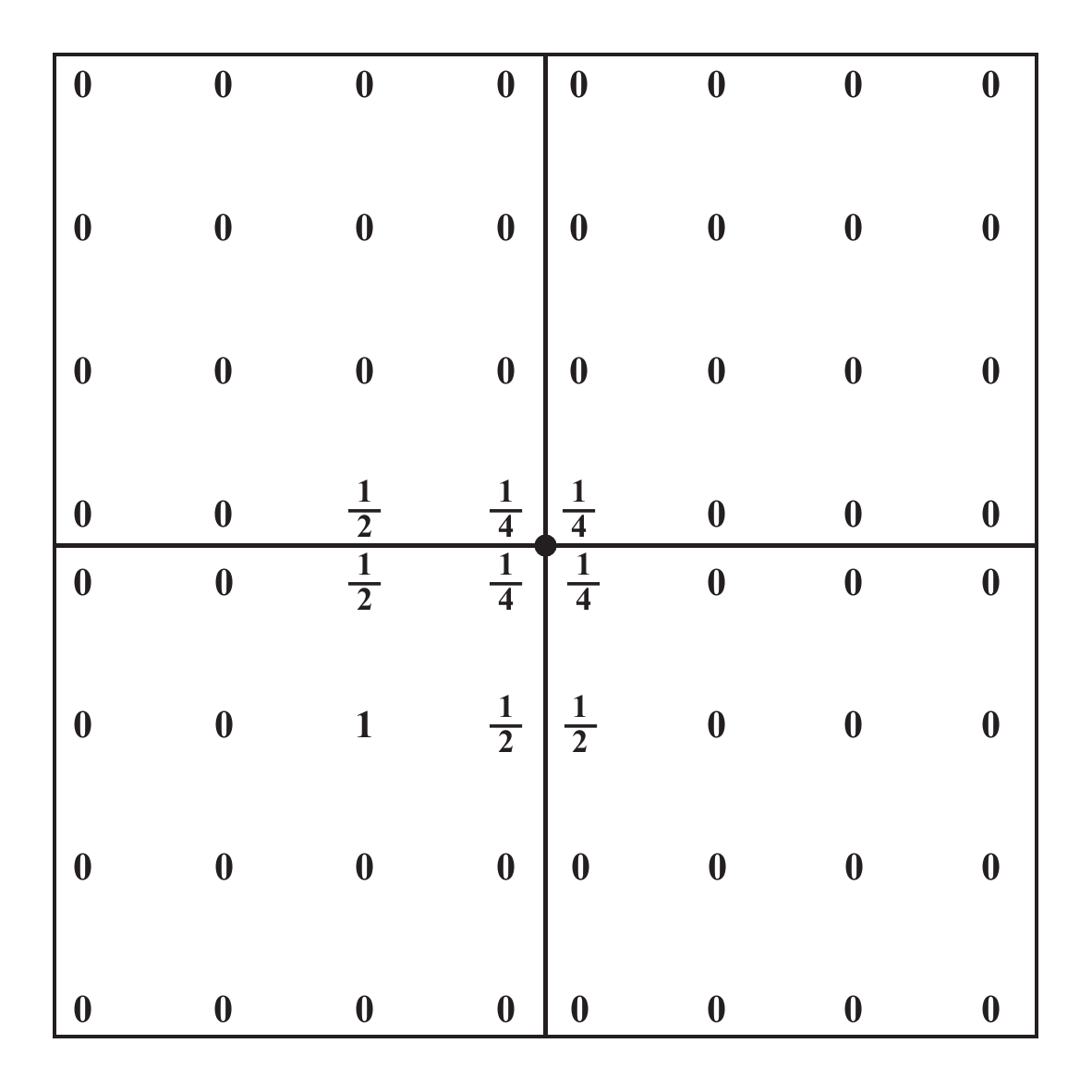}} 
	\subfigure[ B\'{e}zier ordinates on the refined mesh  ]{
		\includegraphics[width = .3\textwidth]{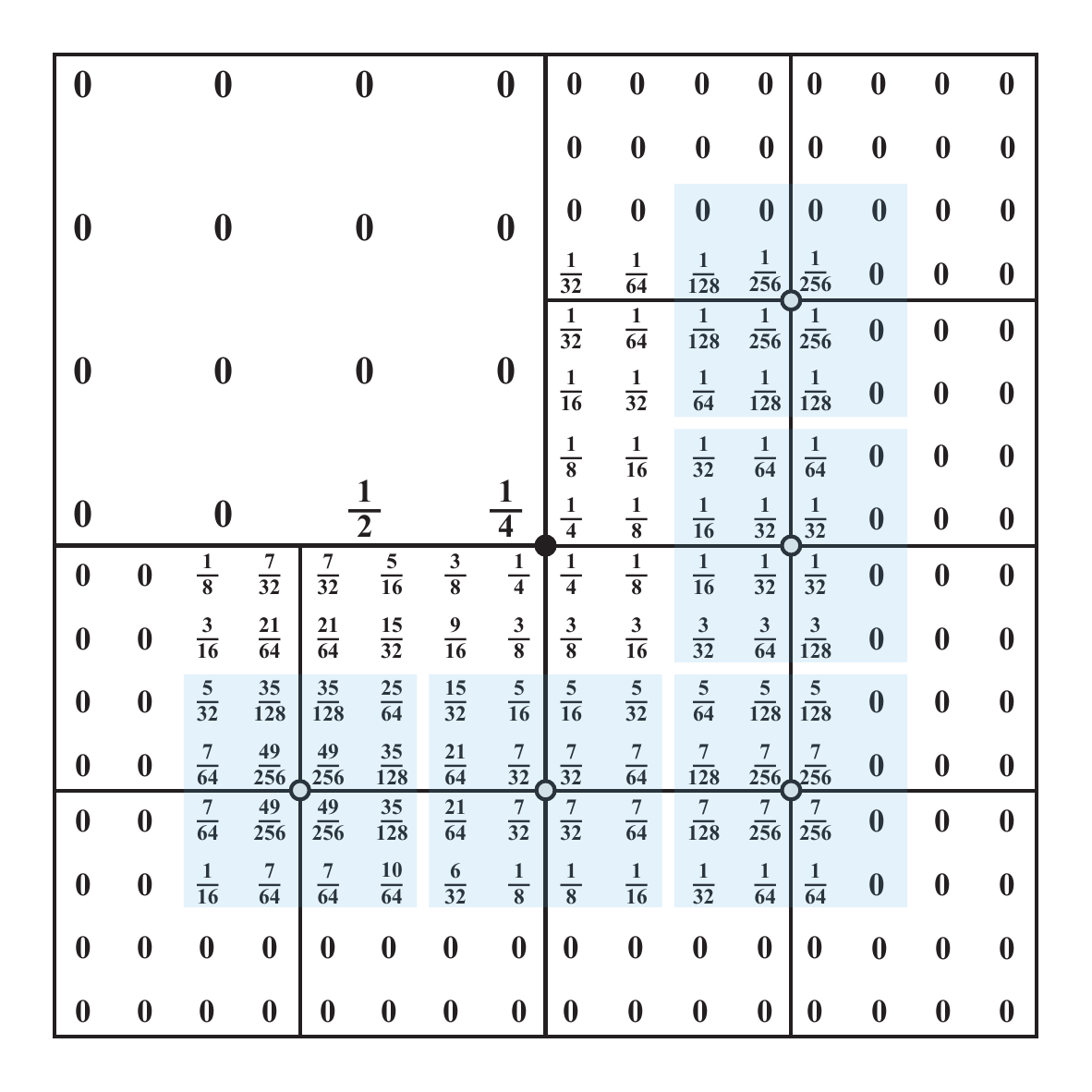}
	} \subfigure[Set B\'{e}zier ordinates associated with new basis
	vertices to be zero.]{
		\includegraphics[width = .3\textwidth]{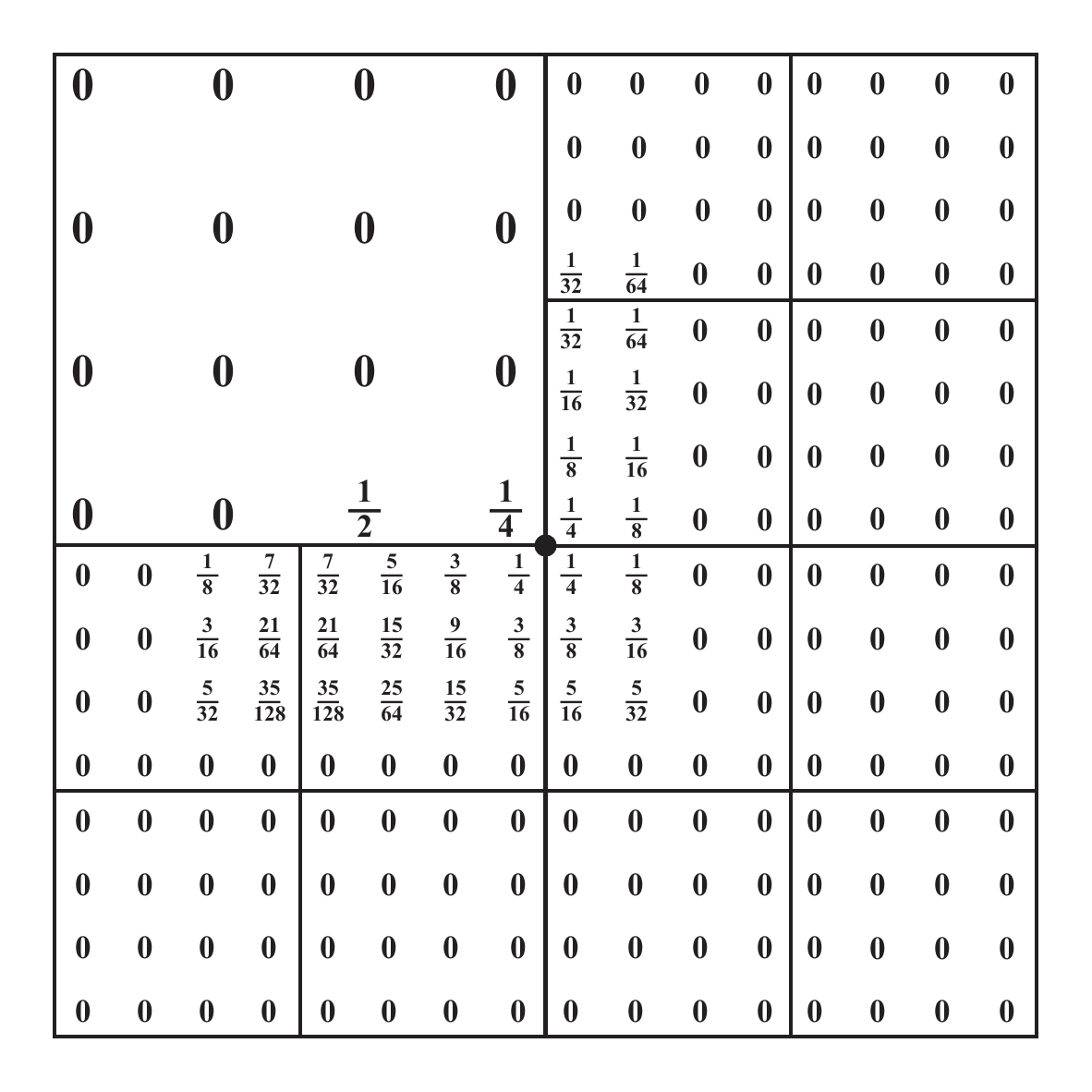}
	} \caption{Illustration of the construction of PHT-spline basis functions. (a) and (b) show the changes in basis functions as the mesh is refined. The last three plots take one of the four basis functions shown in (a) as an example to show how to construct basis functions on the refined mesh by modifying the B\'{e}zier ordinates. The new basis vertices are marked with black circles. The B\'{e}zier ordinates in the shaded rectangles are set to be zeros.  }\label{F.phtbasis}
\end{figure*}

Now we describe the details of the first step. PHT-spline basis functions are represented in B\'{e}zier form by specifying its $16$ B\'{e}zier ordinates in every cell within the support of the basis function, referring to Fig.~\ref{F.phtbasis}(c). From $\mt_k$ to $\mt_{k+1}$, some cells are refined and some new basis vertices appear. The B\'{e}zier coordinates on these refined cells are subdivided into four parts, and each part is associated with a cell corner vertex as shown in Fig.~\ref{F.phtbasis}(d). Then all the basis functions $\{b^k_j\}_{j=1}^{d_k}$ at level $k$ are modified to $\{\overline{b}^k_j\}_{j=1}^{d_k}$ in the following fashion: for each $j$, reset all the associated B\'{e}zier ordinates with the new basis vertices to zero.

The basis functions constructed as described are not only nonnegative and linearly independent, but they also form a partition of unity while possessing localized support. In this paper, the topology optimization of shell structures is conducted in the isogeometric analysis framework based on PHT-splines. The mid-surface, density function, and displacement field of a shell structure are all represented using PHT-splines.

\subsection{Shell structures}
\begin{figure}[!htbp]
	\centering
	\includegraphics[width=0.45\textwidth]{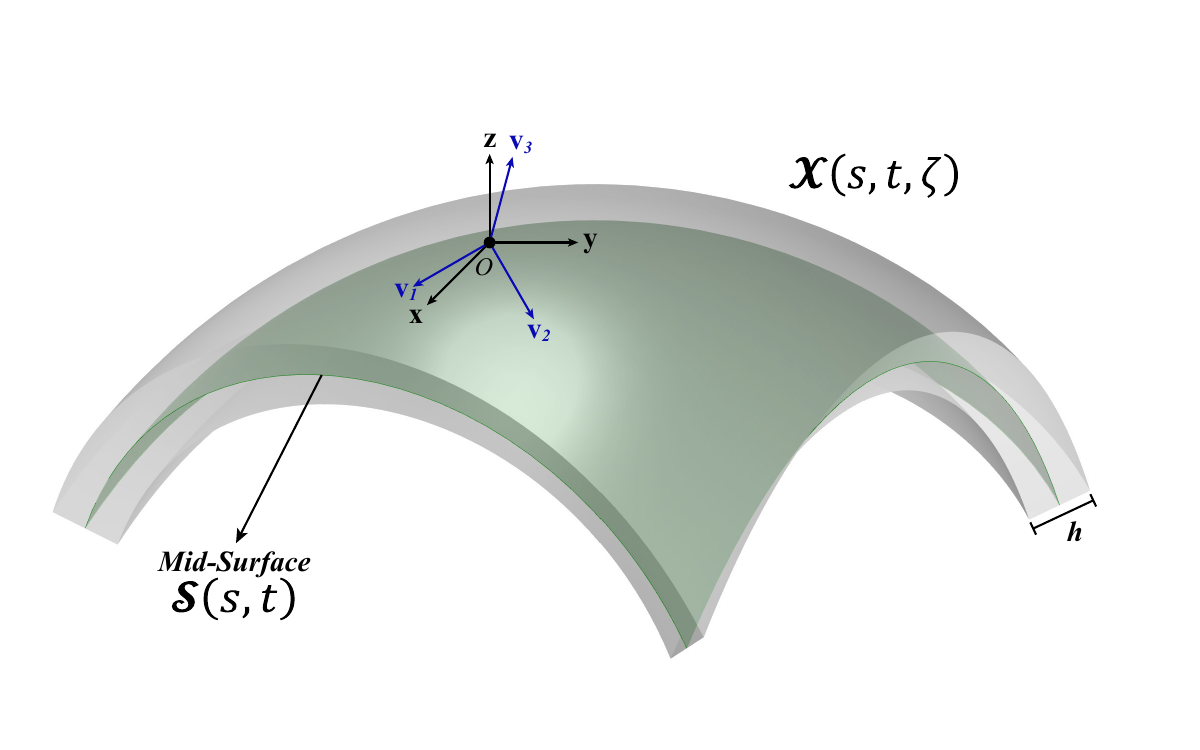}
	\caption{Shell structure, where $\mathcal{S}(s,t) $ is the mid-surface,  $\mathcal{X}(s,t,\zeta) $ is the shell structure, $h$ is the thickness, $\{\bf{v}_1, \bf{v}_2, \bf{v}_3\}$ is the local coordinate system at the point $O$. Especially, $\bf{v}_1$, $\bf{v}_2$ are the two tangent vectors of the mid-surface and $\bf{v}_3$ is the normal at the point $O$. }
	\label{fig:shellstructure}
\end{figure}

A shell structure $\rm \mathcal{X}$ is described by its mid-surface $\mathcal{S}(s,t) $ together with its thickness (as shown in Fig.~\ref{fig:shellstructure}). 
The mid-surface $\mathcal{S}(s,t)$ is represented as a PHT-spline surface on a hierarchical T-mesh $\mt$ with
\begin{equation}
	\label{eq:midsurf}
	\mathcal{S}(s,t) = \sum_{i=1}^n N_{i}(s,t) \mathbf{P}_{i},~~(s,t)\in\Omega_0=[0,1]\times[0,1],
\end{equation}
where $N_i$ is the $i$-th PHT-spline basis function defined over $\mt$, $\mathbf{P}_i\in \mathbb{R}^3$ is the $i$-th control point and $n$ is the number of basis functions. In order to express the overall shell structure in PHT-spline form, a local system on the mid-surface is required. For any point $O$ on the mid-surface, the local system $\{\bf{v}_1,\bf{v}_2,\bf{v}_3\}$ is constructed as in Fig.~\ref{fig:shellstructure}. The unit normal vector at  point $O$ is defined by the cross product of the first order derivatives from $\mathcal{S}$, that is,
\begin{equation}\label{eq:v3}
	{\bf v}_3(s,t)  =
	\frac{\mathcal{S}_{s}(s,t)\times \mathcal{S}_{t}(s,t)}{\|\mathcal{S}_{s}(s,t)\times \mathcal{S}_{t}(s,t)\|},
\end{equation}
Where $\mathcal{S}_{s}=\partial\mathcal{S}/\partial s$ and $\mathcal{S}_{t}=\partial\mathcal{S}/\partial t$. The shell analysis in this paper is based on the most widely used Reissner-Mintlin theory \cite{shell-theory}, which assumes zero normal stress. Based on this assumption, the stress along the direction of ${\bf{v}}_3$ is assumed to be zero.  Once $\bf{v}_3$ is determined, ${\bf{v}}_1$ can be chosen as the outer product of ${\bf{v}}_3$ and one of the three axis vectors in the global system $\{x,y,z\}$. Subsequently, ${\bf{v}}_2$ can be taken as the outer product of ${\bf{v}}_1$ and ${\bf{v}}_3$ \cite{kang2015isogeometric}.

\begin{figure}[!htbp]
	\centering
	\includegraphics[width=0.6\textwidth]{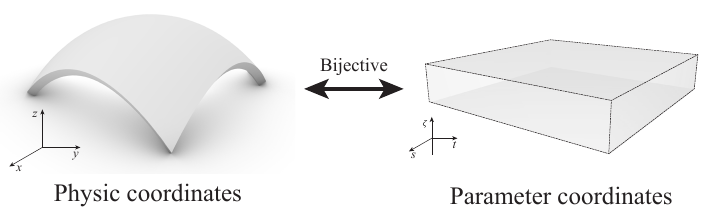}
	\caption{There is a bijective mapping between physical coordinates $(x,y,z)$ and parameter coordinates $(s,t,\zeta)$.}
	\label{F.bijective}
\end{figure}

Now the shell structure is formulated as
\begin{equation}
	\label{eq:shellrepre}
	\mathcal{X}(s,t,\zeta) = \mathcal{S}(s,t) + \zeta~{\bf v}_3(s,t)~\sum_{i=1}^n \frac{h}{2}N_{i}(s,t)= \left(x(s,t,\zeta),y(s,t,\zeta),z(s,t,\zeta)\right)^T,
\end{equation}
where $\mathcal{S}(s,t)$ is the mid-surface defined by Eq.\eqref{eq:midsurf}, $h$ is a constant and defines the thickness of a shell structure, $\zeta \in [-1, 1]$ is the parameter along the thickness direction. This mapping gives a bijective mapping from the parameter coordinates $(s,t,\zeta)\in \Omega_0\times[-1,1]$ to the physical coordinates $(x,y,z)$, as shown in Fig.~\ref{F.bijective}.

The displacement and density field of shell structures are established over hierarchical T-meshes derived from the adaptive refinement of the underlying mesh of the mid-surface. Various levels of these hierarchical T-meshes are denoted as $\mt_k$. Field variables defined on the physical domain (i.e., shell structure, mid-surface) can be retrieved from the PHT-splines defined over the parameter domain via the utilization of the inverse mappings $\mathcal{X}^{-1}$ or $\mathcal{S}^{-1}$.

The displacement $\mathcal{U}(s,t)$ is described by the PHT-spline defined over the hierarchical T-mesh $\mt_k$,  
\begin{equation}
	\label{E.displacement}
	\mathcal{U}(s,t, \zeta) =
	\sum_{i=1}^{n_k}N_{i}^{(k)}(s,t)
	\begin{pmatrix}
		u_i\\v_i\\w_i
	\end{pmatrix}
	+ \zeta ~\frac{h}{2}~{\boldsymbol{\mu}}(s,t)\sum_{i=1}^{n_k}N_i^{(k)}(s,t)
	\begin{pmatrix}
		\alpha_i \\ \beta_i
	\end{pmatrix},
\end{equation}
where $N_i^{(k)}(s,t)$ is the basis function defined over $\mt_k$, $n_k$ is the number of PHT basis on $\mt_k$, $(u_i, v_i, w_i)$ is the translations along the directions of $x-$axis, $y-$axis and $z-$axis, respectively. The variable ${\boldsymbol{\mu}}(s,t)$ denotes the two local coordinate axis ${\boldsymbol{\mu}}(s,t)= (-{\bf v}_2(s,t), {\bf v}_1(s,t))$. The variable $(\alpha_i, \beta_i)$ denotes the rotation angle at the projection point around the local coordinates axis $\bf v_1$ and $\bf v_2$, respectively. In this setting, the displacement $\mathcal{U}(s,t,\zeta)$ is determined by $5$ degrees of freedom per basis function, consisting of three translation variables and two rotation variables with respect to two tangential directions $\bf v_1$ and $\bf v_2$. The displacement of the shell structure is thus defined by $\mathcal{U}\circ \mathcal{X}^{-1}(x,y,z)$.

The objective of topology optimization of shell structures is to minimize compliance while maintaining the volume of filled material below a predefined fraction. Conventional density-based topology optimization aims to identify the optimal material distribution layout within a designated domain by utilizing the density of each finite element as the design variable. In this paper, density is expressed as continuous PHT-spline functions defined on hierarchical T-meshes, specifically 
\begin{equation}
	\label{E.density}
	\mathscr{\rho}(s,t)=\sum_{i=1}^{n_k} N_i^{(k)}(s,t)~\rho_i,
\end{equation}
where $\rho_i\in \mathbb{R}$, $i=1,\cdots,n_k$ are density control coefficients whose values are in the range of $[0,1]$. The density function $\mathscr{\rho}$ shares the same hierarchical T-meshes as the displacement field $\mathcal{U}$.
The density function of the mid-surface is thus defined by $\mathscr{\rho}\circ \mathcal{S}^{-1}(x,y,z)$. The final density $\Tilde{\rho}$ of a shell structure is determined by
\begin{equation}
	\Tilde{\rho}=  \left\{
	\begin{aligned}
		&\text{void}, ~~~ &\rho(s,t) < 0.5,\\
		&\text{solid}, ~~~ &\rho(s,t) \geq 0.5.
	\end{aligned}
	\right.
\end{equation}

\section{Adaptive isogeometric topology optimization of shell structures}
\label{sec:method}
\begin{figure}[!htbp]
	\centering
	\includegraphics[width=\textwidth]{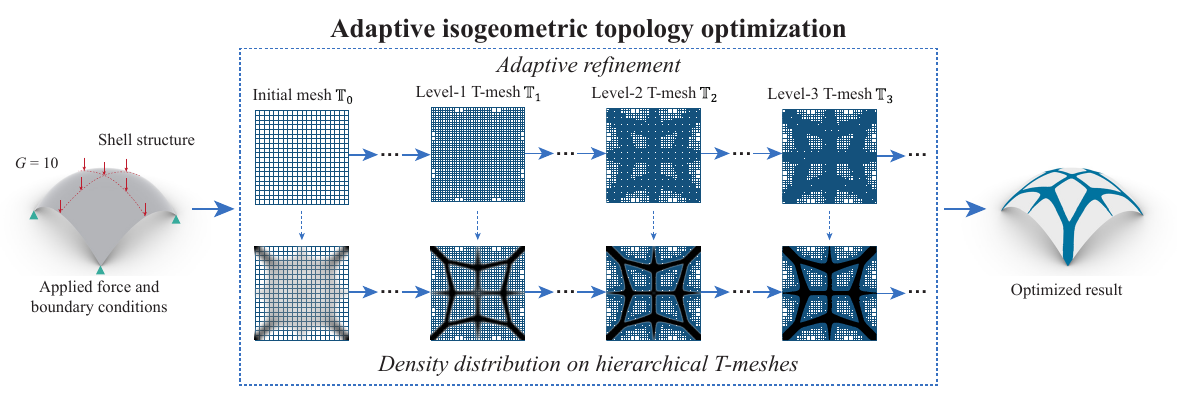}
	\caption{The pipeline of PHT-AITO. }
	\label{fig:pipeline}
\end{figure}

Fig.~\ref{fig:pipeline} shows the PHT-AITO pipeline, which consists of two key components: shell topology optimization and adaptive refinement. The shell topology optimization process starts with a coarse mesh and iteratively performs ITO on nested hierarchical T-mesh sequences. The initial design variable values on each level are inherited from the optimized results on the previous mesh. The control coefficients of the density function serve as design variables of the shell topology optimization. Mesh elements undergo refinement when their density falls within the range of $[\rho_{l}, \rho_{u}]$. The element density is defined by the value of the density function at the element center. In the following, the PHT-AITO is introduced in detail. 

% Mesh elements undergo refinement when their density falls within the range of $[0.1, 0.9]$. The element density is defined by the value of the density function at the element center.

% \textcolor{red}{1. Refinement is performed when the elemental density lies within the specified range of $[0.1, 0.9]$. 2. The refinement condition is not clear. Besides, refinement is performed on T-meshes(parametric space) but the density is a property of physical elements on the shell. There is an ambiguity. Check all the element related variables and mesh definition in the context.}

\subsection{Discretization of TO problem based on PHT-splines}
For a given shell structure $\mathcal{X}$ represented by PHT-splines, the displacement field and density function are established according to Eq.\eqref{E.displacement} and Eq.\eqref{E.density}, respectively. The material distribution in the shell structure is determined by solving the TO problem~(\ref{eq:topopt}) utilizing PHT-splines discretization on hierarchical T-meshes. These meshes are created through the adaptive refinement of the initial mesh $\mt_0$. The TO formulation on the $k$-th level hierarchical T-mesh $\mt_k$ is defined as
\begin{equation}
	\label{eq:topopt}
	\begin{aligned}
		\min_{\{\rho_1^{(k)},\cdots,\rho_{n_k}^{(k)}\}} C &= (\textit{\textbf{U}}^{(k)})^\mathsf{T} \textit{\textbf{K}}(\rho)^{(k)} \textit{\textbf{U}}^{(k)},\\
		s.t. ~~  \textit{\textbf{K}}( \mathscr{\rho})^{(k)}\textit{\textbf{U}}^{(k)}&= \textit{\textbf{F}}^{(k)},\\
		V &\leq V^*,\\
		0 \leq &\rho_i^{(k)} \leq 1,
	\end{aligned} 
\end{equation}
where $\{\rho_{1}^{(k)}, ..., \rho_{n_k}^{(k)}\}$ are the control coefficients of the density function $\rho$ defined over $\mt_k$,  $n_k$ is the number of PHT basis functions over $\mt_k$, $C$ is the compliance of the optimized shell structure, $\textit{\textbf{F}}^{(k)}$ is the force vector corresponding to the shell load,  $\textit{\textbf{U}}^{(k)} = \{(u_{i}^{(k)}, v_{i}^{(k)}, w_{i}^{(k)}, \alpha_{i}^{(k)}, \beta_{i}^{(k)}),i=1,...,n_k\}$ is the control coefficient vector of the displacement $\mathcal{U}(s, t, \zeta)$ defined over $\mt_k$, $\textit{\textbf{K}}(\rho)^{(k)}$ is the stiffness matrix on level $k$, $V^*$ is the predefined volume fraction, $V$ is the volume of the shell structure. Here, we use the superscript $(k)$ to express the variables related to $k$-th level mesh.

The stiffness matrix $\textit{\textbf{K}}_{e}^{(k)}$ on the element $e=[s_{i_0},s_{i_1}]\times[t_{j_0},t_{j_1}]\times[-1,1]$ is approximated by
\begin{equation}
	\label{eq.stiffmat}
	\textit{\textbf{K}}_{e}^{(k)} := \int_{e} (\textit{\textbf{B}}^{(k)})^\mathsf{T} \textit{\textbf{D}} \textit{\textbf{B}}^{(k)}|\textit{\textbf{J}}|dsdtd\zeta   
	\approx E(\rho_e)\textit{\textbf{K}}_{0, e}^{(k)},
\end{equation}
where the Jacobian matrix $\textit{\textbf{J}}$ and the strain-displacement matrix $\textit{\textbf{B}}^{(k)}$ are presented in \textbf{Appendix}, the symbol $\rho_e:=\rho(e)$ denotes the material density at the center of element $e$, $\textit{\textbf{D}}$ is the strain-stress matrix related to Young's modulus $E(\rho)$, and $E(\rho)$ is  defined by the SIMP model,
\begin{equation}
	\label{E.Yongmod}
	E(\rho) = E_{min} + \rho^{p}(E_0 - E_{min}),
\end{equation}
with $E_{min}$ being the minimal material stiffness used to avoid singularity in the stiffness matrix, $E_0$ denotes the exact Young's modulus for the isotropic material, $p$ is the penalization power 
and set to be $p=5$ in this paper. The element stiffness matrix $\textit{\textbf{K}}_{0, e}^{(k)}$ is defined by
\begin{equation}
	\label{eq.stiffmate0}
	\textit{\textbf{K}}_{0, e}^{(k)} := \int_{e} (\textit{\textbf{B}}^{(k)})^\mathsf{T} \textit{\textbf{D}}_0 \textit{\textbf{B}}^{(k)} |\textit{\textbf{J}}|dsdtd\zeta,
\end{equation}
where the $\textit{\textbf{D}}_0$ is strain-stress matrix related to the unit Young's modulus

Note that in the computation of $\textit{\textbf{K}}_{e}^{(k)}$, we simplify the density in an element by the value of the density function $\rho$ at the element center. Throughout each iteration, the stiffness matrix $\textit{\textbf{K}}_{0, e}^{(k)}$ for each element remains fixed, which can be pre-computed in advance. The density-related stiffness matrix $\textit{\textbf{K}}_{e}^{(k)}$ is recalculated with updated density function $\rho$. 
The overall stiffness matrix $\textbf{\textit K}(\rho)^{(k)}$ on $k$-th level is assembled element-wise from $\textbf{\textit K}_e^{(k)}$, which is numerically computed by an integral across the parametric domain $\Omega_0 \times [-1, 1]$. 

The compliance $C$ is computed by
\begin{equation}\label{eq.compliance}
	C = \sum_{e=1}^{n_e}(\textit{\textbf{U}}_e^{(k)})^\mathsf{T}\textit{\textbf{F}}_e^{(k)} =\sum_{i=1}^{n_e}(\textit{\textbf{U}}_e^{(k)})^\mathsf{T} \textit{\textbf{K}}_e^{(k)} \textit{\textbf{U}}_e^{(k)},
\end{equation}
where $n_e$ is the number of elements in the mesh $\mt_k$, $\textit{\textbf{U}}_e^{(k)}$ is the vector of control points for the displacement $\mathcal{U}(s,t,\zeta)$ corresponding to the non-vanishing basis functions on the element $e$, and $\textit{\textbf{F}}_e^{(k)}$ is the force on element $e$. The optimization problem presented in \eqref{eq:topopt} is solved using the Method of Moving Asymptotes (MMA) approach~\cite{mma}.

\subsection{Adaptive smoothing sensitivity}
The sensitivity analysis is performed by taking the gradients of the compliance with respect to the design variables. Here, the control coefficients of the density function \eqref{E.density} are regarded as design variables. The sensitivity filter employed in the classical TO method is defined as the weighted average of the sensitivities within the sphere of a specified radius. However, such a strategy is not suitable for adaptively refined meshes, where element sizes are inconsistent. In \cite{gupta2022adaptive}, it is found that for adaptive meshes, classical filters cause drastic shape changes in subsequent iterations, weakening the optimal topology. In this paper, the sensitivity filter is defined by the weighted average of sensitivities within the support of the target design variable. Here, the support of a design variable refers to the support of the corresponding basis function.

\begin{figure}[!htbp]
	\centering
	\includegraphics[width=0.8\textwidth]{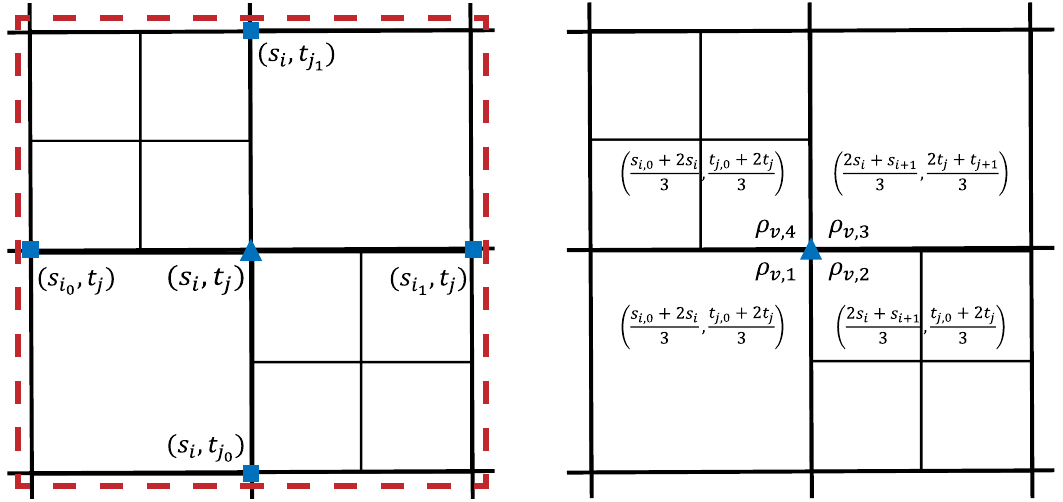}
	\caption{ Left: the support mesh corresponding to the basis vertex marked by a triangle is marked by dashed lines. Right: the parametric coordinates of the four associated density control coefficients.\label{fig:coord}  }
\end{figure}

For the design variable $\rho_i$, the design variables $\rho_j$ within its support are assigned weights based on the Euclidean distance from $\rho_i$, and their sensitivities are filtered as follows:
\begin{equation}
	\label{eq.filter}
	\frac{\partial C}{\partial \rho_i} = \frac{1}{\rho_i\sum_{j=1}^{N}{w_j}}\sum_{j=1}^{N}w_j\rho_j\frac{\partial C}{\partial \rho_j},
\end{equation}
where $C$ is the compliance on the current mesh, $N$ is the number of the design variables within the support of $\rho_i$, and $w_j$ is the weight which is calculated by the Shepard function \cite{shepard1968two}
\begin{equation*}
	w_j = (1-r)^{6}_{+}(35r^2 + 18r + 3),
\end{equation*}
with $r$ being the distance between $\rho_i$ and $\rho_j$.

The distance between $\rho_i$ and $\rho_j$ is defined by the distance between their parametric coordinates. Note that one basis vertex is associated with four design variables (see Sec. \ref{sec.pht}). To distinguish these four design variables, we introduce the support mesh to define parametric coordinates. For a basis vertex $v=(s_i,t_j)$, the \emph{support mesh} is defined as the minimal $2\times 2$ tensor product mesh whose central vertex is $v$ and edges are composed of the edges in meshes. Denote the domain occupied by the support mesh by $[s_{i_0},s_{i_1}]\times[t_{j_0},t_{j_1}]$. Then, the parametric coordinates of the four associated design variables, denoted by $\{\rho_{v,1}, \rho_{v,2}, \rho_{v,3}, \rho_{v,4}\}$, are defined as follows,
\begin{equation}
	\label{eq:coord}
	\begin{aligned}
		\rho_{v,1} &= \left(\frac{s_{i_0} + 2s_{i} }{3}, \frac{t_{j_0} + 2t_{j} }{3}\right),&~
		\rho_{v,2} &= \left(\frac{2s_{i} + s_{i_1}}{3}, \frac{t_{j_0} + 2t_{j}}{3}\right),~\\
		\rho_{v,3} &= \left(\frac{2s_{i} + s_{i_1}}{3}, \frac{2t_{j} + t_{j_1}}{3}\right),&~
		\rho_{v,4} &= \left(\frac{s_{i_0} + 2s_{i}}{3}, \frac{2t_{j} + t_{j_1}}{3}\right).
	\end{aligned}
\end{equation} 
Refer to Fig.~\ref{fig:coord} as an example. The support mesh of $v$ is actually the minimal tensor product mesh that contains the support of the basis functions associated with $v$. Due to the hierarchical structure of the T-mesh, the support mesh for each basis vertex is easily obtained.

\begin{figure}[!htbp]
	\centering
	\includegraphics[width=0.4\textwidth]{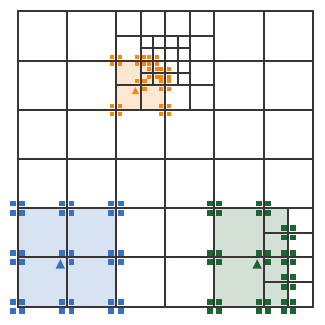}
	\caption{Determination of neighboring design variables. The target design variables are marked by blue, green and yellow triangles. The neighboring design variables are marked as squares with the same color as the corresponding target design variable. 
		\label{fig:filter}}
\end{figure}

The neighboring design variables $\rho_j$ in the sensitivity filter \eqref{eq.filter} are taken as those associated with basis vertices within the support of $\rho_i$. Fig.~\ref{fig:filter} illustrates how neighboring design variables can be determined, where the target $\rho_i$ is marked with a triangle, the corresponding local support region is shaded by color, and neighboring design variables are marked by squares.

Note that the filter radius of \eqref{eq.filter} is automatically determined by the support region of the target design variable, eliminating the requirement for a predefined filter radius. In the meanwhile, the filter radius adaptively changes as the mesh is refined.

\subsection{Adaptive mesh refinement}
The mesh refinement aims to produce finer elements along the material boundaries and coarser elements within the solid and void regions, eliminating redundant elements and degrees of freedom caused by global refinement. Therefore, the gray elements, located along the interface of solid and void regions, are identified for refinement. An element $e$ is marked as a gray element if 
\begin{equation}\label{eq:record}
	\rho_{l}< \rho_e < \rho_{u},
\end{equation}   
where $\rho_e$ is the element density defined by the value of the density function $\rho$ at the element center, $\rho_{l}$ and $\rho_{u}$ are the lower and upper bounds, and they are set as $\rho_{l}=0.1$ and $\rho_u=0.9$ in this paper. 

\begin{figure}[!htbp]
	\centering
	\includegraphics[width=0.9\textwidth]{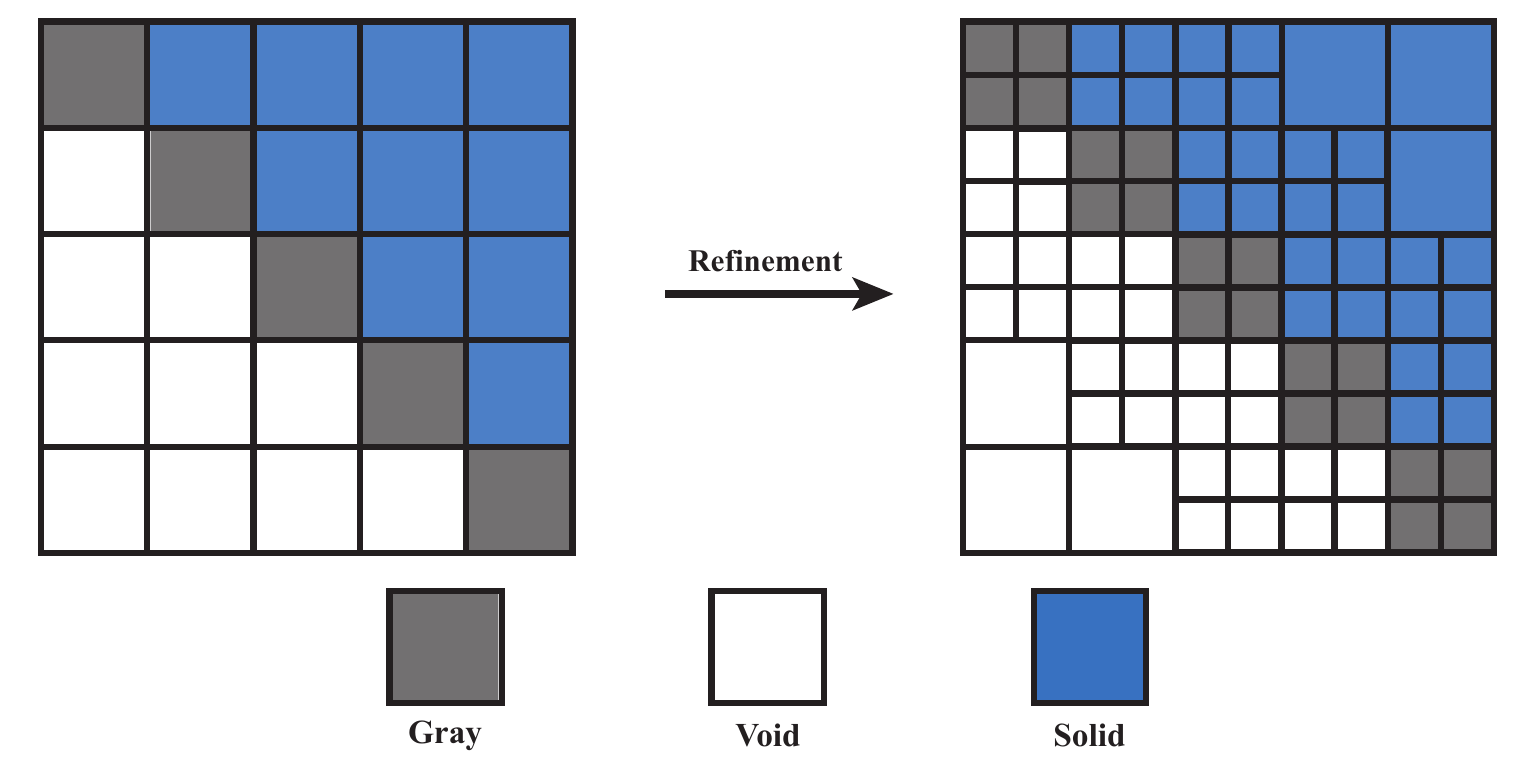}
	\caption{The refinement strategy. The gray elements and solid elements are shaded with gray and blue, respectively.}
	\label{fig:refine}
\end{figure}

Each gray element is divided into four equal sub-elements originating from the center point. In practice, only refining gray elements fails to comprehensively capture the material boundaries. Therefore, when an element undergoes refinement, those neighboring elements that share either an edge or a vertex with the element under consideration, are also refined.
Referring to Fig.~\ref{fig:refine}, the elements on the diagonal are gray elements, and the gray elements together with their neighboring elements are refined.

The algorithmic implementation is illustrated in \textbf{Algorithm }\ref{alg:refine}.
The algorithm runs whenever the maximum value of the iterative variation in element density is greater than the specified tolerance. In order to stabilize the topology structure throughout the optimization process and avoid high computational costs, we avoid mesh refinement at every iteration step. Through some experiments, We found that mesh refinement is initiated when the maximal change of the density coefficients in three consecutive iterations is less than the threshold and the relative change in compliance is less than a predefined tolerance. This strategy provides a good balance between solution accuracy and computational cost. The Relative Change of Compliance (RCC)  is defined by
\begin{equation}
	RCC = \frac{|C_{new} - C_{old}|}{C_{old}},
\end{equation}
where $C_{new}$ is the compliance calculated for the current iteration at the current level mesh, and $C_{old}$ is the compliance at the first iteration on the current level mesh. A larger $RCC$ indicates that further mesh refinement would help improve the optimization results. Therefore, the $RCC$ together with the iterative variation of density is used in the iterative process to restrict redundant mesh refinement, thereby reducing computational costs.

\begin{algorithm}
	\caption{Algorithmic implementation of PHT-AITO}
	\label{alg:refine}
	\hspace*{0.02in} {\bf Input:} 
	$\mathcal{S}$: The mid-surface of shell model represented by PHT-splines; \\
	~~~~~~~~~~~~$\mt_0$: The initial mesh for optimization;\\
	\hspace*{0.02in} {\bf Output:} 
	$\mathbf{\rho}^*$: The vector of control coefficients of density $\rho$.\\
	\begin{algorithmic}
		\STATE Initialize the vector of control coefficients of density $\mathbf{\rho}_{old} = \mathbf{1}$, ~$tol_c=0.01$, ~$tol_{ref}=0.025$, ~$tol_{rcc}=0.15$, ~$count=0$, ~$\rho_{l}=0.1$, ~$\rho_{u}=0.9$, ~$k=0$,  ~$C_{old}=1\times10^{-5}$, $ch=1.0$;
		\WHILE {$ch >= tol_c$}
		\STATE Compute $\textbf{\textit K}^{(k)}$ and $\textbf{\textit F}^{(k)}$ on the hierarchical T-mesh $\mt_{k}$ with $\rho_{old}$ based on Eq.\eqref{eq.stiffmat} and \eqref{eq.load}, respectively;\\
		\STATE Solve $\textbf{\textit U}^{(k)}$ from the linear system $\textbf{\textit K}^{(k)}\textbf{\textit U}^{(k)}=\textbf{\textit F}^{(k)}$;\\
		\STATE Compute the compliance $C_{new}$ based on Eq.\eqref{eq.compliance};\\
		\IF{Compliance $C_{new}$ is first evaluated on $\mt_k~(k\geq1)$}
		\STATE $C_{old} = C_{new}$;\\
		\ENDIF
		\STATE Apply filter for adaptive smoothing sensitivity based on Eq.\eqref{eq.filter};\\
		\STATE Solve TO problem Eq.\eqref{eq:topopt} using MMA and obtain $\rho_{new}$;\\
		\STATE $ch = \max(|\rho_{new} - \rho_{old}|)$;\\
		\STATE $\rho_{old} = \rho_{new}$,  ~$\rho^* = \rho_{new}$;\\
		\STATE $RCC = |C_{new} - C_{old}|/C_{old}$;\\
		\IF{$ch<tol_{ref}$}
		\STATE $count = count + 1$;\\
		\ENDIF
		\IF{$count==3$ \&\& $RCC>tol_{rcc}$}
		\FOR{e $\in \mt_{k}$ }
		\IF{$\rho_{l} < \rho_e < \rho_{u}$}
		\STATE Mark element $e$ as a gray element;\\
		\ENDIF
		\ENDFOR
		\STATE Refine the marked elements and their neighboring elements to obtain the refined mesh $\mt_{k+1}$;
		\STATE Refresh the representation of the density $\rho$ defined by control coefficient vector $\rho_{old}$ on the refined T-mesh $\mt_{k+1}$, according to Eq.\eqref{e.rhoinherit};\\
		\STATE $count = 0$;\\
		\STATE $k = k + 1$;\\
		\ENDIF
		\ENDWHILE
		\STATE Output $\rho^*$;
	\end{algorithmic}
\end{algorithm}

\subsection{Density inheritance}
Within the PHT-AITO framework, density is represented by PHT spline functions and undergoes two distinct updates. The first involves adjusting the control coefficients using MMA, and the second involves increasing the resolution of density representation through mesh refinement. Note that providing initial values of the control coefficient is a requirement in MMA. However, an issue arises during mesh refinement due to the variation in the PHT-spline space. This leads to varying dimensions of the control coefficient across different levels of meshes, making the previously optimized coefficient unsuitable as initial values for MMA on the refined mesh. The straightforward method of initializing all values on the refined mesh as ones necessitates repeating the TO problem on that mesh from scratch. This approach wastes the optimized results that have been obtained on the previous mesh and disregards the inherent connection between the density function across various mesh levels. To address this issue economically and stably, we represent the density defined on $\mt_k$ based on the refined mesh $\mt_{k+1}$. This can be easily achieved since the PHT-spline spaces defined on hierarchical T-meshes are nested, allowing for the exact expression of the density on $\mt_k$ within $\mt_{k+1}$.

Let $\rho^{(k)}$ denote the density defined over $\mt_k$. Our objective is to calculate the control coefficients $\rho_i^{(k+1)}$ for the representations of $\rho^{(k)}$ on $\mt_{k+1}$ such that it satisfies:
\begin{equation}
	\rho^{(k)} = \sum_{i=1}^{n_{k+1}} \rho_i^{(k+1)} N_i^{(k+1)}(s,t).
\end{equation}
The geometric information introduced in \cite{PHT} is used for the conversion. The geometric information of a function $f(s,t)$ is defined by 
\begin{equation}
	\mathcal{L}(f(s,t)) = \left[f, ~\frac{\partial f}{\partial s}, ~\frac{\partial f}{\partial t}, ~\frac{\partial^2 f}{\partial s \partial t} \right].
\end{equation}
For a basis vertex $v=(s_0, t_0)$ in $\mt_{k+1}$, all the PHT-spline basis functions $N_i^{k+1}(s,t)$ vanish at $v$ except the four basis functions associated with $v$, based on this property, we have 
\begin{equation}
	\mathcal{L}(\rho^{(k)}(s_0, t_0)) = \sum_{i\in \mathcal{I}}\rho_i^{k+1} \mathcal{L}(N_i^{k+1}(s_0, t_0)) = \mathcal{P}\cdot\mathcal{G},
\end{equation}
where $\mathcal{I}=\{i_0,i_1,i_2,i_3\}$ is the index set that contains the indices of the four basis functions associated with $v$, $\mathcal{P}=(\rho^{k+1}_{i_0}, \rho^{k+1}_{i_1}, \rho^{k+1}_{i_2}, \rho^{k+1}_{i_3})$ and $\mathcal{G}$ is a $4\times 4$ matrix with the element being the geometric information of basis functions. Then the new density coefficients are computed as
\begin{equation}
	\label{e.rhoinherit}
	\mathcal{P} = \mathcal{L}(\rho^{(k)}(s,t))\cdot \mathcal{G}^{-1}.
\end{equation}
Note that $\mathcal{G}^{-1}$ can be computed explicitly \cite{PHT}, contributing to efficient density inheritance.

\begin{figure}[!htbp]
	\centering
	\subfigure[Density distribution on the coarse mesh]{\includegraphics[width=0.25\textwidth]{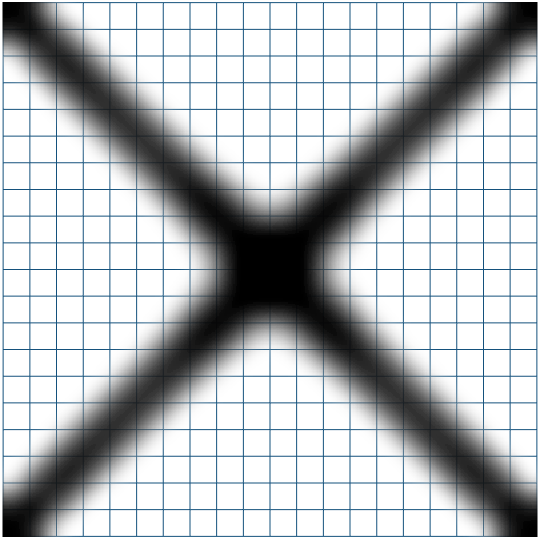}}
	\hspace{3cm}
	\subfigure[The inherited density distribution on the refined mesh]{\includegraphics[width=0.25\textwidth]{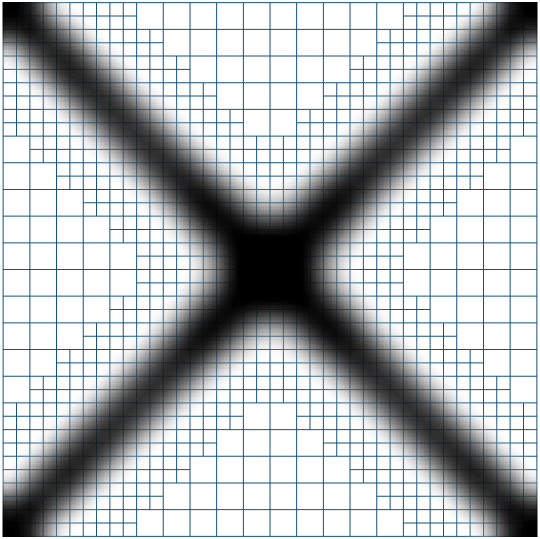}}
	\caption{The density on the refined mesh inherits the density on the coarse mesh, that is, the same density is displayed on two meshes with different resolutions.  }
	\label{fig:inherit}
\end{figure}

Fig. (\ref{fig:inherit}) illustrates the concept of density inheritance. In Fig. (\ref{fig:inherit})(a), we depict the density defined over a uniform mesh. This density is subsequently inherited on the refined mesh, as depicted in Fig. (\ref{fig:inherit})(b). This inheritance ensures that the density on the uniform mesh is exactly preserved in the representation on the refined mesh. Consequently, the topology optimization problem on the refined mesh utilizes these control coefficients from this representation as the initial values for the design variables associated with the refined mesh.

\section{Numerical experiments and discussion}
\label{sec:res}
In this section, numerical experiments on various shell structures are presented to demonstrate the superiority of PHT-AITO. Firstly, we comprehensively compare PHT-AITO with the isogeometric topology optimization method based on B-splines (B-ITO) to demonstrate the computational efficiency and robustness of PHT-AITO.  Subsequently, we provide a detailed presentation of the PHT-AITO process. Furthermore, a comparison is presented among the initial mesh, the final T-mesh, and the hierarchical T-meshes, confirming the hierarchical T-meshes as the preferable choice for shell topology optimization. Additionally, a comparison between density inheritance and resetting density is conducted to underscore the necessity of density inheritance.  Finally, we discuss the influence of two parameters: the initial mesh and the refinement tolerance ($tol_{ref}$ in \textbf{Algorithm} \ref{alg:refine}). These parameters are demonstrated to exert significant influence on the optimized results.

In the following numerical experiments, we use the settings $E_0 = 2100$, $E_{min} = 1\times10^{-5}$, the Poisson's ratio $v = 0.3$ and the refinement threshold $tol_{ref} = 0.15$. The shell thickness is set as $h = 5$. Mesh elements are refined when their density value is in the range of $[0.1,0.9]$. The optimization process terminates when either the maximum change in density coefficients $ch$ is less than $0.01$ or after 300 iterations. The PHT-AITO implementation is implemented in C++ and executed on a laptop with an Intel i9 CPU @ 2.50GHz and 16GB of RAM.

\subsection{Comparison}\label{sec:comparison}
With PHT-AITO, the density distribution on shell structures can be robustly optimized in high-resolution spaces, addressing the problem of excessive freedom caused by full high-resolution meshes by solving topology optimization problems on adaptively refined meshes. To demonstrate this advantage, we conduct a comparison between the proposed PHT-AITO and B-ITO on ten distinct shell structures in terms of degrees of freedom (DoF) and computational time. With B-ITO, the ITO problem is solved by cubic B-splines defined on tensor product meshes without mesh refinement. The proposed algorithm is tested and compared on ten shell structures with different loads and boundary conditions shown in Fig.~\ref{fig:cases}. Details of the shell structure are provided in \textbf{Appendix}.

Noting that the basis functions in PHT-spline space are cubic $C^1$ continuous splines, for fairness of the comparison, the B-ITO framework also uses the cubic B-splines with $C^1$ continuity as the representation of displacement and density. We are considering two types of B-splines. The first one is defined over the tensor product mesh that is obtained by globally refining the initial mesh used in PHT-AITO. This type is called B-splines-1. The second type is defined over the tensor product mesh with degrees of freedom close to the final hierarchical T-mesh obtained from PHT-AITO, called B-splines-2.

\begin{figure}[!htbp]
	\centering
	\includegraphics[width=\textwidth]{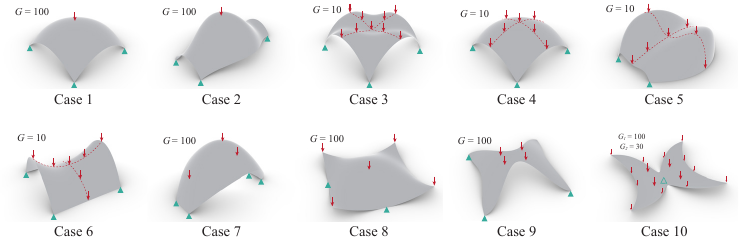}
	\caption{Ten shell structures with boundary conditions and loads.}
	\label{fig:cases}
\end{figure}

\begin{table}[!htbp]
	\centering
	\caption{The number of elements and basis functions, the volume (Vol), the compliance (C), and the computational time of the optimized structures obtained by PHT-AITO and B-ITO. }
	\vspace{0.2cm}
	\label{tab:comparison}
	\renewcommand\arraystretch{1.4}
	\resizebox{!}{.45\textheight}{
		\begin{tabular}{ccccccccc}
			\hline
			\multicolumn{2}{c}{\multirow{2}{*}{Case}} & \multirow{2}{*}{Elements} & \multicolumn{2}{c}{Basis} & \multirow{2}{*}{Vol} & \multirow{2}{*}{Compliance} & \multicolumn{2}{c}{Time(s)} \\ \cline{4-5} \cline{8-9} 
			\multicolumn{2}{c}{}                      &                           & Number     & Ratio      &                      &                             & Number         & Ratio       \\ \hline
			\multirow{3}{*}{1}          & PHT         & 3268                      & 12212      &           & 0.299                & 10.768                      & 184.265        &            \\\cline{2-9}
			& B-1         & $6400(80\times 80)$       & 26244     & 2.149      & 0.299                & 10.734                      & 429.332        & 2.329       \\
			& B-2         & $3136(56\times 56)$       & 12996      & 1.064      & 0.299                & 11.048                      & 232.281        & 1.261       \\ \hline\hline
			\multirow{3}{*}{2}          & PHT         & 3088                      & 11716      &           & 0.299                & 10.198                      & 127.589        &            \\\cline{2-9}
			& B-1         & $6400(80\times 80)$       & 26244     & 2.240      & 0.300                & 10.062                      & 675.355        & 5.293       \\
			& B-2         & $2916(54\times 54)$       & 12100      & 1.033      & 0.300                & 10.503                      & 179.450        & 1.406       \\ \hline\hline
			\multirow{3}{*}{3}          & PHT         & 12568                     & 47292     &           & 0.299                & 176.333                     & 1320.095       &            \\\cline{2-9}
			& B-1         & $25600(160\times 160)$    & 103684     & 2.192      & 0.299                & 178.966                     & 4925.644       & 3.731       \\
			& B-2         & $11664(108\times 108)$    & 47524     & 1.005      & 0.302                & 193.059                     & 1460.771       & 1.106       \\ \hline\hline
			\multirow{3}{*}{4}          & PHT         & 12568                     & 47316     &           & 0.300                & 154.087                     & 1777.319       &           \\\cline{2-9}
			& B-1         & $25600(160\times 160)$    & 103684     & 2.191      & 0.299                & 183.687                     & 6834.035       & 3.845       \\
			& B-2         & $11664(108\times 108)$    & 47524     & 1.004      & 0.303                & 175.271                     & 1948.057       & 1.096       \\ \hline\hline
			\multirow{3}{*}{5}          & PHT         & 12700                     & 47740     &          & 0.298                & 977.166                     & 1581.799       &           \\\cline{2-9}
			& B-1         & $25600(160\times 160)$    & 103684     & 2.172      & 0.300                & 977.219                     & 6918.273       & 4.373       \\
			& B-2         & $12100(110\times 110)$    & 49284     & 1.032      & 0.299                & 1046.018                    & 2260.722      & 1.429       \\ \hline\hline
			\multirow{3}{*}{6}          & PHT         & 12256                     & 46116     &           & 0.297                & 242.896                     & 916.059        &            \\\cline{2-9}
			& B-1         & $25600(160\times 160)$    & 103684     & 2.248      & 0.301                & 262.951                     & 6928.720       & 7.563       \\
			& B-2         & $11664(108\times 108)$    & 47524     & 1.031      & 0.300                & 258.716                     & 1999.325       & 2.182       \\ \hline\hline
			\multirow{3}{*}{7}          & PHT         & 12388                     & 46164     &         & 0.300                & 117.159                     & 783.380        &            \\\cline{2-9}
			& B-1         & $25600(160\times 160)$    & 103684     & 2.246      & 0.298                & 122.838                     & 5975.509       & 7.627       \\
			& B-2         & $11664(108\times 108)$    & 47524     & 1.029      & 0.302                & 138.743                     & 2335.888       & 2.981       \\ \hline\hline
			\multirow{3}{*}{8}          & PHT         & 11596                     & 43636     &           & 0.300                & 610.710                     & 1688.850       &            \\\cline{2-9}
			& B-1         & $25600(160\times 160)$    & 103684     & 2.376      & 0.299                & 596.681                     & 6929.355       & 4.103       \\
			& B-2         & $10816(104\times 104)$    & 44100     & 1.011      & 0.302                & 655.843                     & 2148.908       & 1.272       \\ \hline\hline
			\multirow{3}{*}{9}          & PHT         & 10204                     & 37844      &           & 0.299                & 241.541                     & 992.339        &            \\\cline{2-9}
			& B-1         & $25600(160\times 160)$    & 103684     & 2.739      & 0.300                & 255.112                     & 6901.835       & 6.955       \\
			& B-2         & $9604(98\times 98)$       & 39204      & 1.036     & 0.298                & 245.926                     & 2032.499       & 2.048       \\ \hline\hline
			\multirow{3}{*}{10}         & PHT         & 11164                     & 42684     &           & 0.300                & 1116.171                    & 1180.063       &            \\\cline{2-9}
			& B-1         & $25600(160\times 160)$    & 103684     & 2.429      & 0.298                & 1138.455                    & 6826.129       & 5.784       \\
			& B-2         & $10816(104\times 104)$    & 44100     & 1.033      & 0.300                & 1370.142                    & 2315.048       & 1.962       \\ \hline\hline
	\end{tabular}}
\end{table}

\begin{figure}[!htbp]
	\centering
	\includegraphics[width=0.95\textwidth]{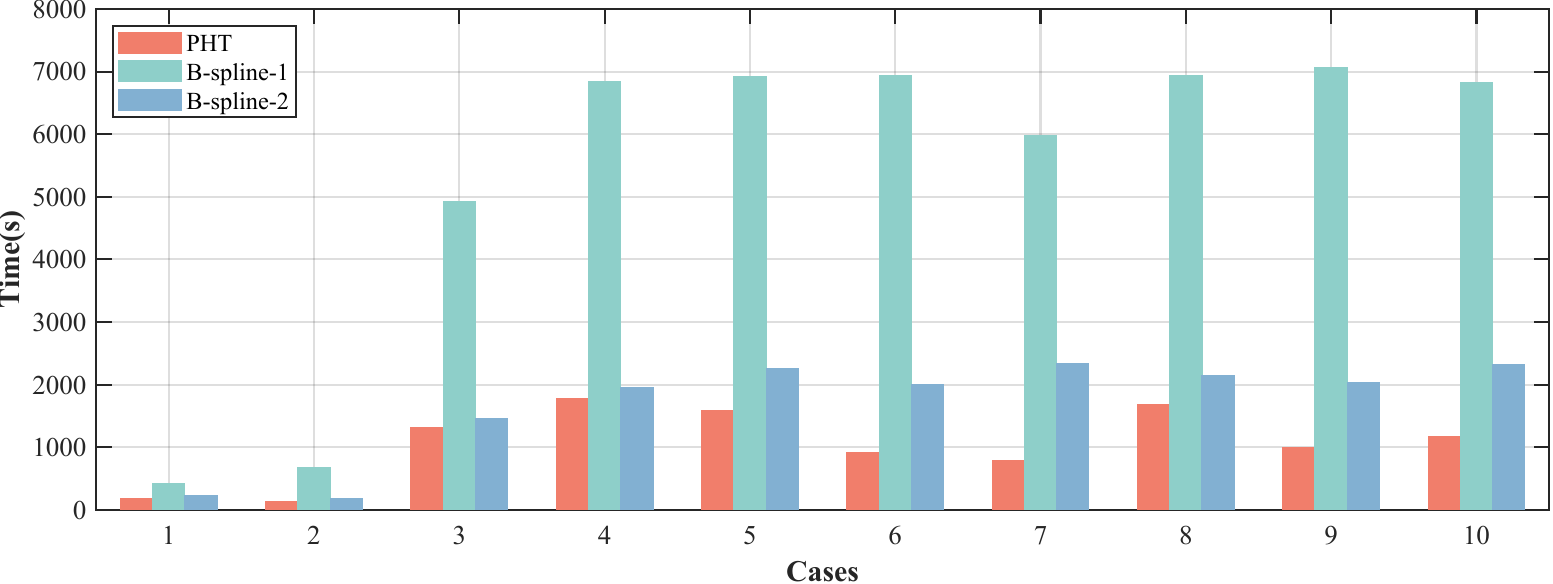}
	\caption{Time-cost bar chart from PHT-AITO, B-ITO based on B-splines-1 and B-splines-2. }
	\label{fig:comparison-time}
\end{figure}

\begin{figure}[!htbp]
	\centering
	\resizebox{!}{.95\textheight}{
		\includegraphics{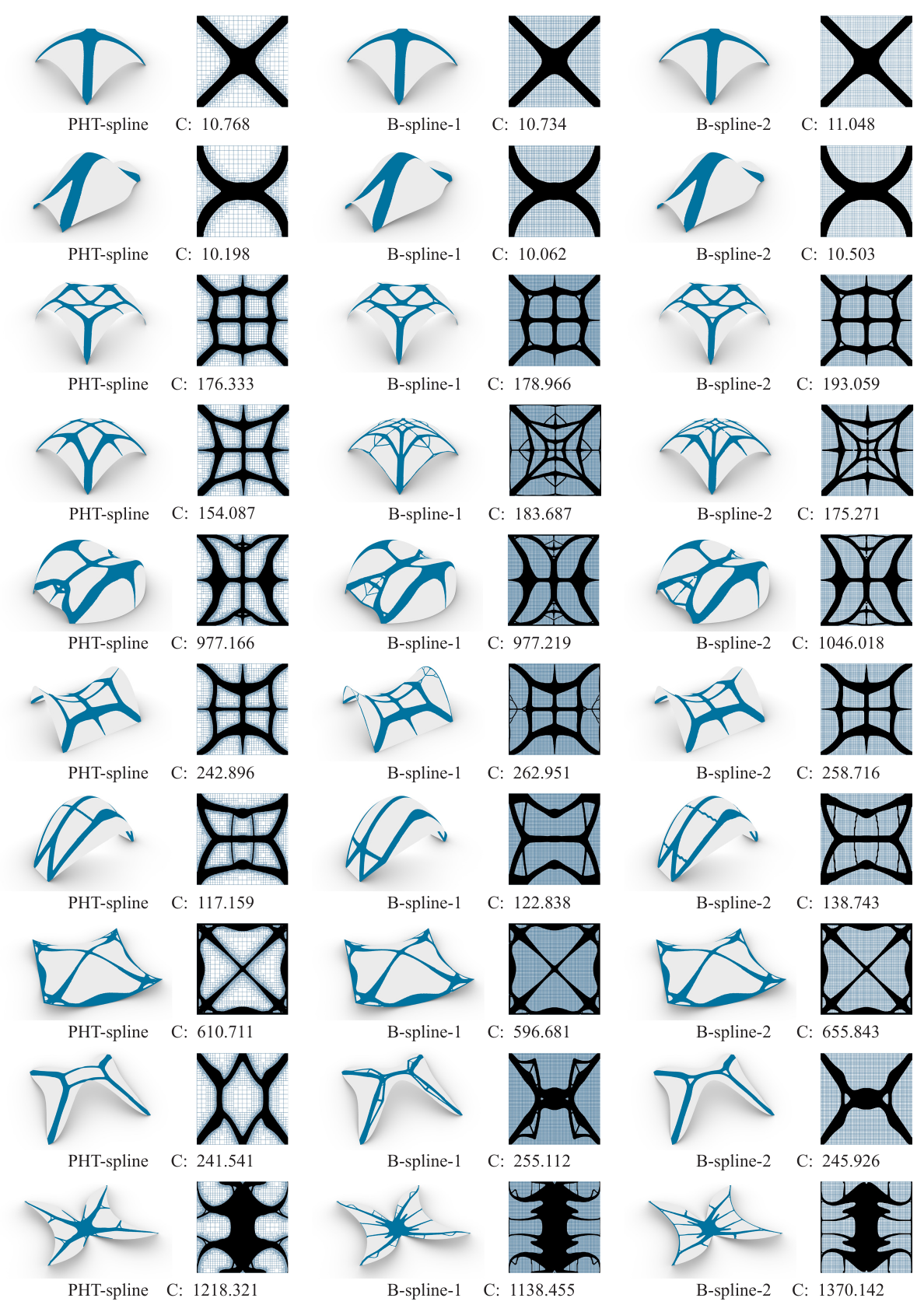}
	}
	\caption{The left, middle, and right two columns show the optimized structures and the density distribution on the corresponding meshes solved by PHT-AITO, B-ITO based on B-splines-1 and B-splines-2, respectively. The symbol $C$ represents the compliance of the optimized structures.}
	\label{fig:comparison-top}
\end{figure}

The statistics of the comparisons are listed in Table~\ref{tab:comparison}. 
For fairness, the volume and compliance of the optimized structures obtained from PHT-AITO and B-ITO are evaluated on a tensor product mesh with a resolution of $200\times 200$. The ratio in the Basis column and Time column is equal to the quotient of the number of basis functions and time consumption used in B-ITO and that used in PHT-AITO, respectively.
The time cost of PHT-AITO and B-ITO is also shown in the form of a bar chart in Fig.~\ref{fig:comparison-time}. We observe that PHT-AITO produces optimized structures at the same resolution as B-ITO but with only half the degrees of freedom needed by B-ITO. Moreover, time efficiency improves by a factor of four on average. Hence, PHT-AITO effectively reduces the redundant DoFs associated with high-resolution representation, achieving greater efficiency.

Moreover, we also compare PHT-AITO with the B-ITO based on B-splines-2. The B-splines-2 are designed to have similar DOFs to the final PHT-splines solved from PHT-AITO, but their underlying meshes are totally different, one is an adaptively refined mesh and the other one is a tensor product mesh. We observe that when the degrees of freedom are nearly identical, PHT-AITO exhibits slightly shorter computation times compared to B-ITO. This suggests that adaptive refinement does not impose additional computational overhead on iterative optimization calculations, instead, it accelerates compliance convergence. We will discuss this advantage further in Section 4.2.

\begin{figure}[!htbp]
	\centering
	\includegraphics[width=0.9\textwidth]{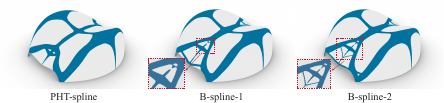}
	\caption{Detailed comparison of the optimized structures for Case 5. The B-ITO obtains optimized structures with fractures, while PHT-AITO can well mitigate these situations.}
	\label{fig:comparison-isolation}
\end{figure}

Fig.~\ref{fig:comparison-top} demonstrates the optimized results obtained by PHT-AITO, and B-ITO based on B-splines-1 and B-splines-2. As can be seen from cases 2-6, the material layouts obtained by PHT-AITO have fewer holes and cracks, while B-ITO prefers to find more intricate material layouts. Case 5 is a representative example of this phenomenon, shown in Fig.\ref{fig:comparison-isolation}, in which the material layouts optimized by B-ITOs have more holes and connections than PHT-AITO. Numerical experiments show that PHT-AITO tends to produce simpler and smoother layouts while B-ITO prefers to produce intricate layouts.

\subsection{The advantage of adaptive optimization }\label{sec:ho}
\begin{figure}[!htbp]
	\centering
	\includegraphics[height=0.45\textheight]{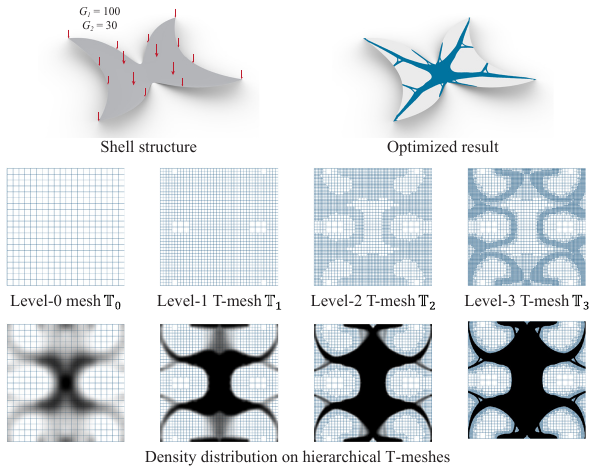}
	\caption{The adaptive topology optimization process of case~10.}
	\label{fig:htop-case10}
\end{figure}

We take case 10 as an example to illustrate the adaptive optimization process. The PHT-AITO framework starts with a $20\times20$ tensor product mesh and iteratively refines the mesh until the stop criterion is satisfied. The iteration process of the mesh refinement and the density distribution is shown in  Fig.~\ref{fig:htop-case10}. There are a lot of gray regions in the density distribution on the level-0 mesh. Through several iterations of mesh adaptive refinement, the gray regions gradually decrease, the material boundaries become much clearer, and a good topology is finally obtained on the level-3 mesh. The refined elements precisely capture the material boundaries well.

A further comparison is conducted to investigate the impact of adaptive refinement on compliance. This analysis encompasses five distinct cases denoted as C1, C2, C3, C4, and C5, specifically,
\begin{itemize}
	\item[-] C1: the ITO problem is solved by the PHT-AITO framework, where the initial mesh is chosen as a $20\times20$ tensor product mesh. Three mesh refinements are conducted in the optimization process.
	\item[-] C2: the ITO problem is solved by the B-ITO framework based on $40\times40$ mesh.
	\item[-] C3: the ITO problem is solved by the B-ITO framework based on $80\times80$ mesh.
	\item[-] C4: the ITO problem is solved by the B-ITO framework based on $104\times104$ mesh, i.e. B-spline-2 in Sec.\ref{sec:comparison}.
	\item[-] C5:the ITO problem is solved by the B-ITO framework based on $160\times160$ mesh, i.e. B-spline-1 in Sec.\ref{sec:comparison}.
\end{itemize}

\begin{figure}[!htbp]
	\centering
	\subfigure[The optimized structures obtained from C1-C5.]{\includegraphics[width=0.75\textwidth]{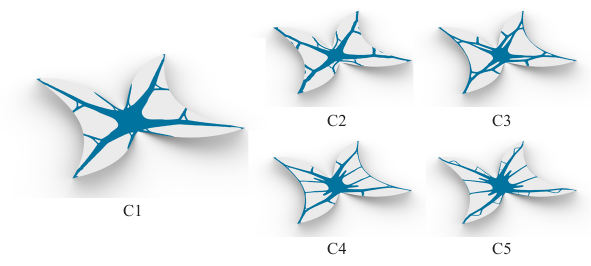}}
	\subfigure[Change curves of compliance over the consumption time]{\includegraphics[width=0.3\textwidth]{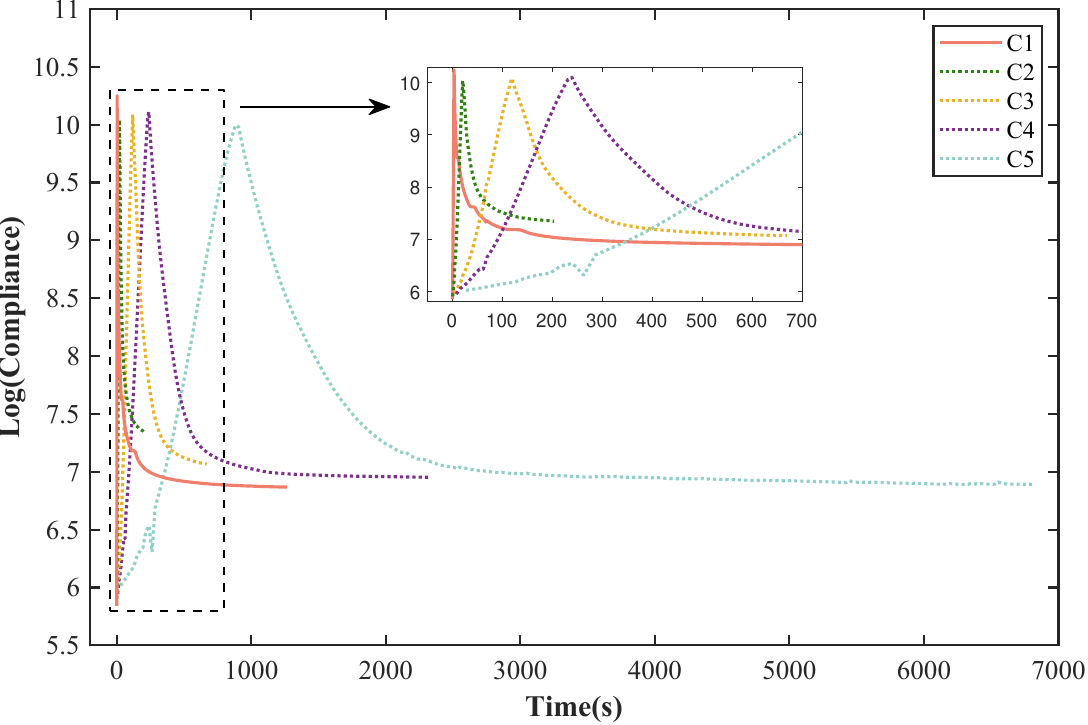}}
	\hspace{1.5cm}
	\subfigure[Variation of the basis number over the iteration progress]{\includegraphics[width=0.3\textwidth]{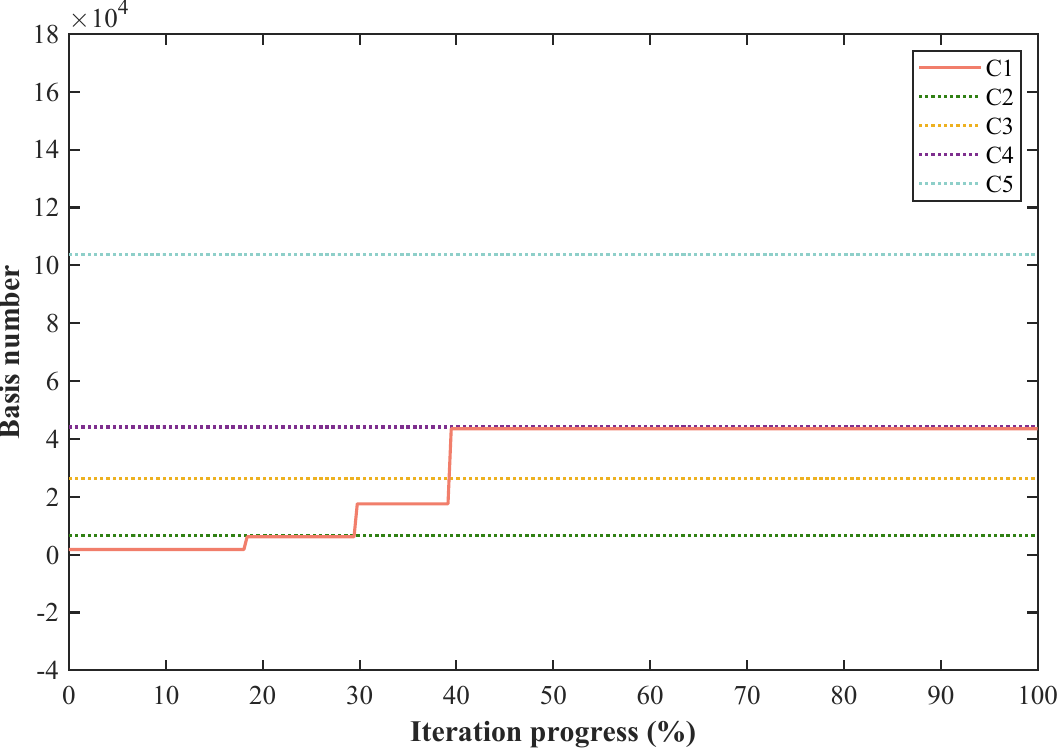}}
	\caption{The optimized structures, compliance, and the number of basis functions in cases C1-C5. }
	\label{fig:optpro-case10}
\end{figure}

The optimized results for five cases are depicted in Fig.~\ref{fig:optpro-case10}(a). The material layouts from C2 have fractures and islands. This may be attributed to the lower-resolution meshes that lack the necessary degrees of freedom to precisely represent the density distribution. As the mesh resolution increases, both C3, C4, and C5 yield more satisfactory material layouts, especially, C5 gives a notably intricate result. Fig.~\ref{fig:optpro-case10}(b) illustrates convergence curves of compliance solved from C1-C5. As the resolution of the tensor product mesh increases (from C2 to C5), compliance progressively converges to lower values. Remarkably, the compliance results from C1 (PHT-AITO) converge faster than those from C4 (B-ITO) with almost the same degrees of freedom. Moreover, the compliance results from C1 converge to almost the same values as C5, but C5 requires twice as many degrees of freedom as C1.

Finally, we investigate the impact of density inheritance on the adaptive optimization process. We compare the compliance obtained by PHT-AITO with density inheritance and density reset. The density reset means the initial values of the density on the refined mesh are all set to be $1$. The change of compliance with density inheritance exhibits a stepwise decline due to the mesh refinement, while the change of compliance has larger oscillations when density reset is used. Additionally, the density reset might result in the refined elements being inconsistent with material boundaries, negating the adaptive advantage of the PHT-AITO. Interestingly, the density reset strategy could produce complex optimized structures with more intricate details, providing a way of constructing diverse optimized results.

\begin{figure}[!htbp]
	\centering
	\includegraphics[width=.8\textwidth]{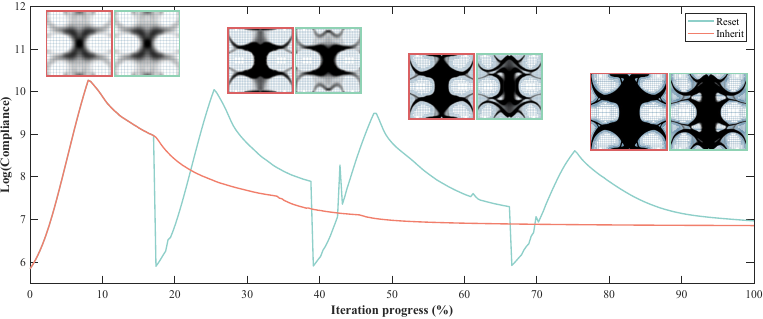}
	\caption{The compliance curves obtained by PHT-AITO using density inheritance and density reset.}
	\label{fig:comparison-DI}
\end{figure}

\subsection{Discussion}

\begin{figure}[!htbp]
	\centering
	\includegraphics[width=0.8\textwidth]{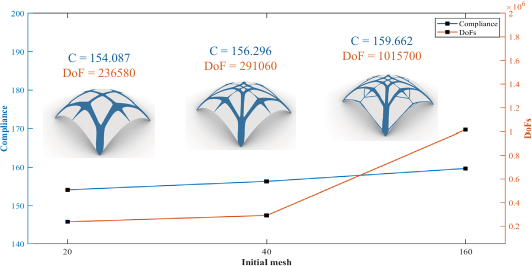}
	\caption{The optimized structures from the PHT-AITO start with initial meshes with different resolutions. }
	\label{fig:initialmesh}
\end{figure}

In the current PHT-AITO framework, lower-resolution tensor product meshes, say $20\times 20$, are chosen as initial meshes, for the sake of stability of topological structures and computational efficiency. Lower-resolution initial meshes tend to produce the profile of a simpler topology structure. By means of density inheritance, it ensures topological consistency and stability during the whole optimization process. Higher-resolution initial meshes tend to produce an intricate material layout and lead to excessive refinement. We take Case~4 as an example to demonstrate the influence of the resolution of initial meshes. Fig.~\ref{fig:initialmesh} shows the optimized shell structures from the PHT-AITO start with initial meshes with different resolutions: $20 \times 20$, $40 \times 40$, and $160 \times 160$. The optimization results exhibit notably distinct topological structures. With the $20 \times 20$ initial mesh, the optimized structure has a simple topology. As the mesh resolution increased, holes and branches occurred in the optimized results.

\section{Conclusion and Future Work}
\label{sec:con}

Shell structures are widely used in architecture and industries. Isogeometric topology optimization enables the creation of innovative material layouts for shell structures. While high-resolution meshes offer more design possibilities, conducting topology optimization on them can be time-consuming. Adaptive Isogeometric Topology Optimization (AITO) addresses the problem of excessive freedom caused by full high-resolution meshes by solving topology optimization problems on adaptively refined meshes. PHT-splines are well-suited for AITO due to their efficient local refinement and advantageous analytical properties. In PHT-AITO, the density function is represented as a PHT-spline function, with control coefficients serving as design variables. These variables are optimized to minimize the compliance of the shell structure under volume constraints. PHT-AITO employs an adaptive mesh for the density function, featuring finer elements at the interface between solid and void material, and coarser elements within the solid and void regions. This approach yields a high-quality density distribution with reduced degrees of freedom. Additionally, density inheritance ensures the stability of the optimized structure and accelerates the iteration process. We validate the effectiveness of PHT-AITO on various shell structures through numerical experiments. These experiments demonstrate that, in contrast to tensor product B-splines, PHT-ATIO offers significantly enhanced computing performance and demonstrates robustness.

We currently provide a viable implementation framework for PHT-AITO, but there is a need to explore specific acceleration and high-performance algorithms for it. One potential way for enhancement is the development of fast matrix assembly techniques and multigrid methods for PHT-splines. Because PHT-splines lack the tensor product structure, existing methods cannot be directly applied, necessitating the development of novel algorithms. Moreover, earlier studies involving ITO-based shell structures simultaneously optimized both shell shape and topology. Hence, a promising direction for further exploration is the extension of PHT-AITO to concurrently address shape and topology optimization.

\section*{Declaration of competing interest}
The authors declare that they have no known competing financial interests or personal relationships that could have appeared to influence the work reported in this paper.

\section*{Acknowledgement}
This work is supported by the Provincial Natural Science Foundation of Anhui (2208085QA01), the Fundamental Research Funds for the Central Universities (WK0010000075), the National Natural Science Foundation of China (No. 61972368,  No. 11801393, No. 12371383) and the
Natural Science Foundation of Jiangsu Province, China (No. BK20180831)

\bibliography{mybibfile}

\begin{thebibliography}{10}
\expandafter\ifx\csname url\endcsname\relax
  \def\url#1{\texttt{#1}}\fi
\expandafter\ifx\csname urlprefix\endcsname\relax\def\urlprefix{URL }\fi
\expandafter\ifx\csname href\endcsname\relax
  \def\href#1#2{#2} \def\path#1{#1}\fi

\bibitem{MICHAELROTTER19983}
J.~{Michael Rotter}, Shell structures: the new {European} standard and current
  research needs, Thin-Walled Structures 31~(1) (1998) 3--23.

\bibitem{bendsoe2003topology}
M.~P. Bendsoe, O.~Sigmund, Topology optimization: theory, methods, and
  applications, Springer Science \& Business Media, 2003.

\bibitem{farshad2013design}
M.~Farshad, Design and analysis of shell structures, Vol.~16, Springer Science
  \& Business Media, 2013.

\bibitem{naghdi1962foundations}
P.~M. Naghdi, Foundations of elastic shell theory, Tech. rep., Institute of
  Engineering Research, University of California Berkeley (1962).

\bibitem{budiansky1968notes}
B.~Budiansky, Notes on nonlinear shell theory, Journal of Applied Mechanics
  35~(2) (1968) 393--401.

\bibitem{ahmad1970analysis}
S.~Ahmad, B.~M. Irons, O.~Zienkiewicz, Analysis of thick and thin shell
  structures by curved finite elements, International journal for numerical
  methods in engineering 2~(3) (1970) 419--451.

\bibitem{hinton1992finite}
E.~Hinton, N.~Rao, J.~Sienz, Finite element structural shape and thickness
  optimization of axisymmetric shells, Engineering computations (1992).

\bibitem{maute1997adaptive}
K.~Maute, E.~Ramm, Adaptive topology optimization of shell structures, AIAA
  journal 35~(11) (1997) 1767--1773.

\bibitem{ansola2002integrated}
R.~Ansola, J.~Canales, J.~A. Tarrago, J.~Rasmussen, An integrated approach for
  shape and topology optimization of shell structures, Computers \& structures
  80~(5-6) (2002) 449--458.

\bibitem{hassani2013simultaneous}
B.~Hassani, S.~M. Tavakkoli, H.~Ghasemnejad, Simultaneous shape and topology
  optimization of shell structures, Structural and Multidisciplinary
  Optimization 48 (2013) 221--233.

\bibitem{Ye2019Topology}
Q.~Ye, Y.~Guo, S.~Chen, N.~Lei, X.~D. Gu, Topology optimization of conformal
  structures on manifolds using extended level set methods {(X-LSM)} and
  conformal geometry theory, Computer Methods in Applied Mechanics and
  Engineering 344 (2019) 164--185.

\bibitem{soto1993modelling}
C.~Soto, A.~Diaz, On the modelling of ribbed plates for shape optimization,
  Structural Optimization 6 (1993) 175--188.

\bibitem{SIMP}
M.~P. Bends{\o}e, O.~Sigmund, Material interpolation schemes in topology
  optimization, Archive of applied mechanics 69 (1999) 635--654.

\bibitem{levelset}
M.~Y. Wang, X.~Wang, D.~Guo, A level set method for structural topology
  optimization, Computer methods in applied mechanics and engineering 192~(1-2)
  (2003) 227--246.

\bibitem{Townsend2019ALS}
S.~Townsend, H.~A. Kim, A level set topology optimization method for the
  buckling of shell structures, Structural and Multidisciplinary Optimization
  60 (2019) 1783 -- 1800.

\bibitem{MMC}
X.~Guo, W.~Zhang, W.~Zhong, Doing topology optimization explicitly and
  geometrically—a new moving morphable components based framework, Journal of
  Applied Mechanics 81~(8) (2014) 081009.

\bibitem{Huo2022topology}
W.~Huo, C.~Liu, Z.~Du, X.~Jiang, Z.~Liu, X.~Guo, Topology optimization on
  complex surfaces based on the moving morphable component method and
  computational conformal mapping, Journal of Applied Mechanics 89~(5) (2022)
  051008.

\bibitem{jiang2023explicit}
X.~Jiang, W.~Zhang, C.~Liu, Z.~Du, X.~Guo, An explicit approach for
  simultaneous shape and topology optimization of shell structures, Applied
  Mathematical Modelling 113 (2023) 613--639.

\bibitem{IGA}
T.~J. Hughes, J.~A. Cottrell, Y.~Bazilevs, Isogeometric analysis: {CAD}, finite
  elements, {NURBS}, exact geometry and mesh refinement, Computer methods in
  applied mechanics and engineering 194~(39-41) (2005) 4135--4195.

\bibitem{cottrell2009isogeometric}
J.~A. Cottrell, T.~J. Hughes, Y.~Bazilevs, Isogeometric analysis: toward
  integration of {CAD} and {FEA}, John Wiley \& Sons, 2009.

\bibitem{seo2010shape}
Y.-D. Seo, H.-J. Kim, S.-K. Youn, Shape optimization and its extension to
  topological design based on isogeometric analysis, International Journal of
  Solids and Structures 47~(11-12) (2010) 1618--1640.

\bibitem{seo2010isogeometric}
Y.-D. Seo, H.-J. Kim, S.-K. Youn, Isogeometric topology optimization using
  trimmed spline surfaces, Computer Methods in Applied Mechanics and
  Engineering 199~(49-52) (2010) 3270--3296.

\bibitem{hassani2012isogeometrical}
B.~Hassani, M.~Khanzadi, S.~M. Tavakkoli, An isogeometrical approach to
  structural topology optimization by optimality criteria, Structural and
  multidisciplinary optimization 45 (2012) 223--233.

\bibitem{gao2019isogeometric}
J.~Gao, L.~Gao, Z.~Luo, P.~Li, Isogeometric topology optimization for continuum
  structures using density distribution function, International Journal for
  Numerical Methods in Engineering 119~(10) (2019) 991--1017.

\bibitem{shojaee2012composition}
S.~Shojaee, M.~Mohamadian, N.~Valizadeh, Composition of isogeometric analysis
  with level set method for structural topology optimization, Int J Optim Civil
  Eng 2~(1) (2012) 47--70.

\bibitem{wang2016isogeometric}
Y.~Wang, D.~J. Benson, Isogeometric analysis for parameterized {LSM-based}
  structural topology optimization, Computational Mechanics 57 (2016) 19--35.

\bibitem{hou2017explicit}
W.~Hou, Y.~Gai, X.~Zhu, X.~Wang, C.~Zhao, L.~Xu, K.~Jiang, P.~Hu, Explicit
  isogeometric topology optimization using moving morphable components,
  Computer Methods in Applied Mechanics and Engineering 326 (2017) 694--712.

\bibitem{zhang2020explicit}
W.~Zhang, D.~Li, P.~Kang, X.~Guo, S.-K. Youn, Explicit topology optimization
  using {IGA}-based moving morphable void {(MMV)} approach, Computer Methods in
  Applied Mechanics and Engineering 360 (2020) 112685.

\bibitem{kang2016isogeometric}
P.~Kang, S.-K. Youn, Isogeometric topology optimization of shell structures
  using trimmed {NURBS} surfaces, Finite Elements in Analysis and Design 120
  (2016) 18--40.

\bibitem{zhang2020stress}
W.~Zhang, S.~Jiang, C.~Liu, D.~Li, P.~Kang, S.-K. Youn, X.~Guo, Stress-related
  topology optimization of shell structures using {IGA/TSA-based} moving
  morphable void {(MMV)} approach, Computer Methods in Applied Mechanics and
  Engineering 366 (2020) 113036.

\bibitem{gao2020comprehensive}
J.~Gao, M.~Xiao, Y.~Zhang, L.~Gao, A comprehensive review of isogeometric
  topology optimization: methods, applications and prospects, Chinese Journal
  of Mechanical Engineering 33~(1) (2020) 1--14.

\bibitem{itobook}
J.~Gao, L.~Gao, M.~Xiao, Isogeometric Topology Optimization: Methods,
  Applications and Implementations, Vol.~7, Springer, 2022.

\bibitem{lrsplines}
X.~Li, F.~Chen, H.~Kang, J.~Deng, A survey on the local refinable splines,
  Science China Mathematics 59 (2016) 617--644.

\bibitem{THBTO}
X.~Xie, S.~Wang, Y.~Wang, N.~Jiang, W.~Zhao, M.~Xu, Truncated hierarchical
  {B-spline–based} topology optimization, Structural and Multidisciplinary
  Optimization 62 (2020) 83--105.

\bibitem{TTO}
G.~Zhao, J.~Yang, W.~Wang, Y.~Zhang, X.~Du, M.~Guo, T-splines based
  isogeometric topology optimization with arbitrarily shaped design domains.,
  Computer Modeling in Engineering \& Sciences 123~(3) (2020) 1033--1059.

\bibitem{gupta2022adaptive}
A.~Gupta, B.~Mamindlapelly, P.~L. Karuthedath, R.~Chowdhury, A.~Chakrabarti,
  Adaptive isogeometric topology optimization using {PHT} splines, Computer
  Methods in Applied Mechanics and Engineering 395 (2022) 114993.

\bibitem{karuthedath2023continuous}
P.~L. Karuthedath, A.~Gupta, B.~Mamindlapelly, R.~Chowdhury, A continuous field
  adaptive mesh refinement algorithm for isogeometric topology optimization
  using {PHT-Splines}, Computer Methods in Applied Mechanics and Engineering
  412 (2023) 116075.

\bibitem{HBiga}
A.-V. Vuong, C.~Giannelli, B.~J{\"u}ttler, B.~Simeon, A hierarchical approach
  to adaptive local refinement in isogeometric analysis, Computer Methods in
  Applied Mechanics and Engineering 200~(49-52) (2011) 3554--3567.

\bibitem{Tsplines}
T.~W. Sederberg, J.~Zheng, A.~Bakenov, A.~Nasri, T-splines and {T-NURCCs}, ACM
  transactions on graphics 22~(3) (2003) 477--484.

\bibitem{AST}
M.~Scott, X.~Li, T.~Sederberg, T.~Hughes, Local refinement of analysis-suitable
  {T}-splines, Computer Methods in Applied Mechanics and Engineering 213-216
  (2012) 206--222.

\bibitem{THB}
C.~Giannelli, B.~Jüttler, H.~Speleers, {THB-splines}: The truncated basis for
  hierarchical splines, Computer Aided Geometric Design 29~(7) (2012) 485--498.

\bibitem{LR}
T.~Dokken, T.~Lyche, K.~F. Pettersen, Polynomial splines over locally refined
  box-partitions, Computer Aided Geometric Design 30~(3) (2013) 331--356.

\bibitem{PHT}
J.~Deng, F.~Chen, X.~Li, C.~Hu, W.~Tong, Z.~Yang, Y.~Feng, Polynomial splines
  over hierarchical {T}-meshes, Graphical models 70~(4) (2008) 76--86.

\bibitem{PHTIGA}
P.~Wang, J.~Xu, J.~Deng, F.~Chen, Adaptive isogeometric analysis using rational
  {PHT}-splines, Computer-Aided Design 43~(11) (2011) 1438--1448.

\bibitem{PHTthinshell}
N.~Nguyen-Thanh, J.~Kiendl, H.~Nguyen-Xuan, R.~Wüchner, K.~Bletzinger,
  Y.~Bazilevs, T.~Rabczuk, Rotation free isogeometric thin shell analysis using
  pht-splines, Computer Methods in Applied Mechanics and Engineering
  200~(47-48) (2011) 3410--3424.

\bibitem{PHTelastic}
N.~Nguyen-Thanh, H.~Nguyen-Xuan, S.~P.~A. Bordas, T.~Rabczuk, Isogeometric
  analysis using polynomial splines over hierarchical {T-meshes} for
  two-dimensional elastic solids, Computer Methods in Applied Mechanics and
  Engineering 200~(21-22) (2011) 1892--1908.

\bibitem{PHTVibrations}
P.~Yu, C.~Anitescu, S.~Tomar, S.~P.~A. Bordas, P.~Kerfriden, Adaptive
  isogeometric analysis for plate vibrations: An efficient approach of local
  refinement based on hierarchical a posteriori error estimation, Computer
  Methods in Applied Mechanics and Engineering 342 (2018) 251--286.

\bibitem{PHTKL}
J.~Videla, F.~Contreras, H.~X. Nguyen, E.~Atroshchenko, Application of
  {PHT-splines} in bending and vibration analysis of cracked {Kirchhoff–Love}
  plates, Computer Methods in Applied Mechanics and Engineering 361 (2020)
  112754.

\bibitem{bib:PHT}
J.~Deng, F.~Chen, X.~Li, C.~Hu, W.~Tong, Z.~Yang, Y.~Feng, Polynomial splines
  over hierarchical {T}-meshes, Graphical Models 74 (2008) 76--86.

\bibitem{shell-theory}
J.~C. Simo, D.~D. Fox, On a stress resultant geometrically exact shell model.
  {Part I}: Formulation and optimal parametrization, Computer Methods in
  Applied Mechanics and Engineering 72~(3) (1989) 267--304.

\bibitem{kang2015isogeometric}
P.~Kang, S.-K. Youn, Isogeometric analysis of topologically complex shell
  structures, Finite Elements in Analysis and Design 99 (2015) 68--81.

\bibitem{mma}
K.~Svanberg, The method of moving asymptotes—a new method for structural
  optimization, International journal for numerical methods in engineering
  24~(2) (1987) 359--373.

\bibitem{shepard1968two}
D.~Shepard, A two-dimensional interpolation function for irregularly-spaced
  data, in: Proceedings of the 1968 23rd ACM national conference, 1968, pp.
  517--524.

\end{thebibliography}

\begin{appendices}
\section*{Supplementary}
\label{sec:supplementary}
\subsection*{Shell structure analysis}
PHT-splines represent the shell structure, displacement field, and density function in PHT-AITO. In the following, we give the deduction of the required stiffness matrix, force, and volume in the problem (\ref{eq:topopt}). For brevity, we have removed the superscript $(k)$ in the following derivation, except for the basis function.

The strain can be expressed in terms of the derivatives of displacement as follows
\begin{equation}
	\label{eq.strain-displacement}
	(\epsilon_x,~\epsilon_y,~\epsilon_z,~\epsilon_{x,y},~\epsilon_{y,z},~\epsilon_{z,x})^\mathsf{T} = \textit{\textbf{H}} (\frac{\partial u}{\partial x},~\frac{\partial u}{\partial y},~\frac{\partial u}{\partial z},~\dots,~\frac{\partial w}{\partial z})^\mathsf{T},
\end{equation}
where $(\epsilon_x, \epsilon_y, \epsilon_z, \epsilon_{x,y}, \epsilon_{y,z}, \epsilon_{z,x})^\mathsf{T}$ is the strain along different directions in physical coordinates system $(x, y,z)$, and $\textit{\textbf{H}}$ is $6\times 9$ constant matrix defined as follows
\begin{equation*}
	\textit{\textbf{H}} = 
	\begin{pmatrix}
		1 & 0 & 0 & 0 & 0 & 0 & 0 & 0 & 0\\
		0 & 0 & 0 & 0 & 1 & 0 & 0 & 0 & 0\\
		0 & 0 & 0 & 0 & 0 & 0 & 0 & 0 & 1\\
		0 & 1 & 0 & 1 & 0 & 0 & 0 & 0 & 0\\
		0 & 0 & 0 & 0 & 0 & 1 & 0 & 1 & 0\\
		0 & 0 & 1 & 0 & 0 & 0 & 1 & 0 & 0\\
	\end{pmatrix}.
\end{equation*}

The partial derivatives of the displacement $\mathcal{U}(s,t,\zeta)=(u,v,w)^T$ with respect to physic coordinates $(x, y, z)$ are calculated as follows,
\begin{equation}
	\label{eq.derivUx}
	\begin{pmatrix}
		\frac{\partial u}{\partial x}\\
		\frac{\partial u}{\partial y}\\
		\frac{\partial u}{\partial z}\\
		\vdots\\
		\frac{\partial w}{\partial z}
	\end{pmatrix} =
	\begin{pmatrix}
		\textit{\textbf{J}}^{-1} & \mathbf{0} & \mathbf{0} \\
		\mathbf{0} & \textit{\textbf{J}}^{-1} & \mathbf{0}\\
		\mathbf{0} & \mathbf{0} & \textit{\textbf{J}}^{-1}
	\end{pmatrix}
	\begin{pmatrix}
		\frac{\partial u}{\partial s}\\
		\frac{\partial u}{\partial t}\\
		\frac{\partial u}{\partial \zeta}\\
		\vdots\\
		\frac{\partial w}{\partial \zeta}
	\end{pmatrix} \triangleq 
	\mathbf{\Gamma}     \begin{pmatrix}
		\frac{\partial u}{\partial s}\\
		\frac{\partial u}{\partial t}\\
		\frac{\partial u}{\partial \zeta}\\
		\vdots\\
		\frac{\partial w}{\partial \zeta}
	\end{pmatrix},
\end{equation}
where $\textit{\textbf{J}}$ is the Jacobian matrix and expressed as
\begin{equation}
	\textit{\textbf{J}} = \begin{pmatrix}
		\frac{\partial x}{\partial s} & \frac{\partial y}{\partial s} &\frac{\partial z}{\partial s}\\
		\frac{\partial x}{\partial t} & \frac{\partial y}{\partial t} &\frac{\partial z}{\partial t}\\
		\frac{\partial x}{\partial \zeta} & \frac{\partial y}{\partial \zeta} &\frac{\partial z}{\partial \zeta}
	\end{pmatrix}.
\end{equation}
The partial derivatives of displacement $(u, v, w)$ with respect to $(s,t,\zeta)$ can be written as
\begin{equation}
	\label{eq.derivUs}
	(\frac{\partial u}{\partial s},~
	\frac{\partial u}{\partial t},~
	\frac{\partial u}{\partial \zeta},~
	\cdots,~
	\frac{\partial w}{\partial \zeta})^T=\textit{\textbf{R}}(
	u_i,~v_i,~w_i,~\alpha_i,~\beta_i)^T.
\end{equation}
Here, $\textit{\textbf{R}} = [\textit{\textbf{R}}_1, \dots, \textit{\textbf{R}}_m]$ with $m$ being the number of non-vanishing PHT basis on the element $e$, and
\begin{equation}
	\textit{\textbf{R}}_i = 
	\begin{pmatrix}
		N_{i,s}^{(k)} & 0       & 0       & -\frac{\zeta h}{2}(m_{1,s}N^{(k)}_i + m_1N^{(k)}_{i,s})  & \frac{\zeta h}{2}(l_{1,s}N^{(k)}_i + l_1N^{(k)}_{i,s}) \\
		N^{(k)}_{i,t} & 0       & 0       & -\frac{\zeta h}{2}(m_{1,t}N^{(k)}_i + m_1N^{(k)}_{i,t})  & \frac{\zeta h}{2}(l_{1,t}N^{(k)}_i + l_1N^{(k)}_{i,t}) \\
		0       & 0       & 0       & -\frac{h}{2}m_1N^{(k)}_i  & \frac{h}{2}l_1N^{(k)}_i \\
		0       & N^{(k)}_{i,s} & 0       & -\frac{\zeta h}{2}(m_{2,s}N_i^{(k)} + m_2N^{(k)}_{i,s})  & \frac{\zeta h}{2}(l_{2,s}N^{(k)}_i + l_2N^{(k)}_{i,s}) \\
		0       & N^k_{i,t} & 0       & -\frac{\zeta h}{2}(m_{2,t}N^{(k)}_i + m_2N^{(k)}_{i,t})  & \frac{\zeta h}{2}(l_{2,t}N^{(k)}_i + l_2N^{(k)}_{i,t}) \\
		0       & 0       & 0       & -\frac{h}{2}m_2N^{(k)}_i  & \frac{h}{2}l_2N^{(k)}_i \\
		0       & 0       & N^{(k)}_{i,s} & -\frac{\zeta h}{2}(m_{3,s}N_i^{(k)} + m_3N^{(k)}_{i,s})  & \frac{\zeta h}{2}(l_{3,s}N^{(k)}_i + l_3N^{(k)}_{i,s}) \\
		0       & 0       & N^{(k)}_{i,t} & -\frac{\zeta h}{2}(m_{3,t}N^k_i + m_3N^{(k)}_{i,t})  & \frac{\zeta h}{2}(l_{3,t}N_i^{(k)} + l_3N^{(k)}_{i,t}) \\
		0       & 0       & 0       & -\frac{h}{2}m_3N^{(k)}_i  & \frac{h}{2}l_3N^{(k)}_i \\
	\end{pmatrix},
\end{equation}
where $N_{i,s}^{(k)}$ and $N_{i,t}^{(k)}$ are the derivatives of the PHT basis $N_{i}^{(k)}$ with respect to $s$ and $t$, respectively. $l_{i, s}$ and $l_{i, t}$ are the derivatives of the $l$ with respect to $s$ and $t$, respectively. The symbols $l_{i,s}, l_{i,t}$, $m_{i,s}$ and $m_{i,t}$ have the similar definition. Additionally, $l_i$ and $m_i$ are the elements of local axis ${\bf v}_1$ and ${\bf v}_2$, i.e. ${\bf v}_1=(l_1,l_2,l_3)$ and ${\bf v}_2=(m_1,m_2,m_3)$.

By combining Eqs.\eqref{eq.strain-displacement}-\eqref{eq.derivUs}, the relationship between strain and displacement is obtained as follows,
\begin{equation}       
	(\epsilon_x,~\epsilon_y,~\epsilon_z,~\epsilon_{x,y},~\epsilon_{y,z},~\epsilon_{z,x})^\mathsf{T} = 
	\textit{\textbf{H}}\mathbf{\Gamma}\textit{\textbf{R}}
	(u_i,~v_i,~w_i,~\alpha_i,~ \beta_i)^\mathsf{T}.
\end{equation}
In the Reissner-Mindlin theory, the shell structure is analyzed on the local system $({\bf v}_1, {\bf v}_2, {\bf v}_3)$ where ${\bf v}_i$, $i=1,2,3$ are the axis of the local system and computed based on Eq.\eqref{eq:v3}. Thus the strain $\epsilon_{local}$ in local system is calculated as 
\begin{equation}       
	\epsilon_{local} \triangleq (\epsilon_{v_1},~\epsilon_{v_2},~\epsilon_{v_3},~\epsilon_{v_1, v_2},~\epsilon_{v_2, v_3},~\epsilon_{v_3,v_1})^\mathsf{T} = 
	\textit{\textbf{T}}\textit{\textbf{H}}\mathbf{\Gamma}\textit{\textbf{R}}
	(u_i,~v_i,~w_i,~\alpha_i,~ \beta_i)^\mathsf{T}, 
\end{equation}
where $\textit{\textbf{T}}$ is the transform matrix from the physical coordinates system to the local system. Therefore, the strain-displacement matrix $\textit{\textbf{B}}$ is calculated as 
\begin{equation}
	\label{eq.B2}
	\textit{\textbf{B}}=\textit{\textbf{T}}\textit{\textbf{H}}\mathbf{\Gamma}\textit{\textbf{R}}.
\end{equation}
Then we substitute \eqref{eq.B2} into \eqref{eq.stiffmate0}, the element stiffness matrix can be computed.

The volume $V$ is computed as
\begin{equation}
	\label{eq.volume}
	V=\sum_{e=1}^{n_e}V_e = \sum_{e=1}^{n_e} \int_e \rho|\textit{\textbf{J}}|dsdtd\zeta \approx\sum_{e=1}^{n_e} V_e^0 \rho_e,
\end{equation}
where $\rho_e$ is the material density at the center of element $e$ and $V_e^0$ denotes the solid volume of element $e$, which is defined by 
\begin{equation}
	\label{E.ve0}
	V_e^0=\int_{e} |\textit{\textbf{J}}|ds dt d\zeta.
\end{equation}

The load vector $\textit{\textbf{F}}$ is calculated by
\begin{equation}
	\label{eq.load}
	\textit{\textbf{F}}=\int_{e} \mathbf{\Psi}G |\textit{\textbf{J}}| dsdtd\zeta,
\end{equation}
where $G$ is the force onto the shell structure, $\textit{\textbf{J}}$ is the Jacobian matrix, and  $\mathbf{\Psi}=[\Psi_{1},...,\Psi_{m}]$ contains all the non-vanishing basis functions on the element $e$, and
\begin{equation}
	\label{eq.psi}
	\Psi_{i}=\begin{pmatrix}
		N_i^{(k)}(s,t) & 0 & 0 & -\frac{h}{2}\zeta m_1 N_i^{(k)}(s,t) & \frac{h}{2}\zeta l_1 N_i^{(k)}(s,t)\\
		0&N_i^{(k)}(s,t) & 0 & -\frac{h}{2}\zeta m_2 N_i^{(k)}(s,t) & \frac{h}{2}\zeta l_2 N_i^{(k)}(s,t)\\
		0&0& N_i^{(k)}(s,t) & -\frac{h}{2}\zeta m_3N_i^{(k)}(s,t) & \frac{h}{2} \zeta l_3 N_i^{(k)}(s,t)\\
	\end{pmatrix}.
\end{equation}

\subsection*{Sensitivity Analysis}

In ITO problem, the objective is to minimize the structure compliance $C = \textit{\textbf{U}}^{\mathsf{T}}\textit{\textbf{K}}\textit{\textbf{U}}$. The sensitivity of compliance $C$ with respect to design variable $\rho_i$ can be written as follows,
\begin{equation}
	\label{eq.sens3}
	\frac{\partial C}{\partial \rho_i} = -\textit{\textbf{U}}^{\mathsf{T}}\frac{\partial \textit{\textbf{K}}}{\partial \rho_i}\textit{\textbf{U}}
	= -\sum_{e}\textit{\textbf{U}}_e^{\mathsf{T}}\frac{\partial \textit{\textbf{K}}_e}{\partial \rho_i}\textit{\textbf{U}}_e,    
\end{equation}
where $\textit{\textbf{U}}_e$ is the control point vector of the displacement field on the element, and $\textit{\textbf{K}}_e$ is the stiffness matrix on the element $e$. The sensitivity can be further written as
\begin{equation}\label{eq.sensC}
	\frac{\partial C}{\partial \rho_i} = -
	\sum_{e} p (\textbf{\textit E}_{0}-\textbf{\textit E}_{min})\rho^{p-1}_e \textbf{\textit U}^T_e \textbf{\textit K}_{0,e} \textbf{\textit U}_e N_i(e),
\end{equation}
where $\rho_e$ is the density value at the center of the element $e$, $N_i(e)$ is the value of the PHT-spline basis $N_i(s,t)$ at the element  center, $\textit{\textbf{K}}_{0, e}$ is computed by \eqref{eq.stiffmate0}.

The sensitivity of the volume \eqref{eq.volume} with respect to $\rho_i$ is defined as follows
\begin{equation}\label{eq.senseV}
	\frac{\partial V}{\partial \rho_i} = 
	\sum_{e} \frac{\partial (V_e^0\rho_e)}{\partial \rho_i}
	= \sum_{e} V_e^0 N_i(e),
\end{equation}
where $\rho_e$ is the density value at the center of the element $e$, and $V_e^0$ is defined by \eqref{E.ve0}.

\subsection*{\textbf{3. Examples of shell structures in Sec.\ref{sec:comparison}}}

The details of the ten shell structures shown in Fig.~\ref{fig:cases} are summarized below. In both cases \textbf{1} and \textbf{2}, the load $G=100$ is applied at the center of the shell with parametric coordinates $(0.5,0.5)$. In cases \textbf{3-6}, the load $G=10$ is applied along two cross directions: the parametric direction $s=0.5$ and $t=0.5$. In cases \textbf{7-9}, the load $G=100$ is applied at several points. In case~7, $G$ is applied at five points with parametric coordinates $(0.2, 0.5)$, $(0.5, 0.2)$, $(0.2, 0.8)$, $(0.8, 0.2)$ and $~(0.5, 0.5)$. In case~8, $G$ is applied at four corners and the center, their parameter coordinates are $(0, 0)$, $(0, 1)$, $(1, 0)$, $(1, 1)$, and $(0.5, 0.5)$. In case~9, the parametric coordinates of applied load $G$ are $(0.2, 0.3)$, $(0.8, 0.7)$, $(0.2, 0.7)$, $(0.8, 0.3)$. In case~10, two kinds of load $G_1=100$ and $G_2=30$ are simultaneously applied on shell, where $G_1$ is applied at $(0.25, 0.33)$, $(0.25, 0.67)$, $(0.75, 0.33)$ and $(0.75, 0.67)$, $G_2$ is applied at $(0,0)$, $(0,1)$, $(1,0)$, $(1,1)$, $(0.25,0)$, $(0.75,0)$, $(0.25,1)$, $(0.75,1)$, $(0,0.4)$, $(0,0.6)$, $(1,0.4)$, $(1,0.6)$. 
Boundary conditions are imposed in this way: the midpoints of the four boundaries are fixed in Case 8, the center of the shell is fixed in Case 10,  and the four corner points are fixed in the rest of cases.

\end{appendices}

\end{document}